\documentclass[12pt]{amsart}
\usepackage{amssymb, latexsym, amsthm}
\usepackage{epsf}
\newtheorem{thm}{Theorem}[section]

\newtheorem{cor}[thm]{Corollary}
\newtheorem{defn}[thm]{Definition}
\newtheorem{exmpl}[thm]{Example}
\newtheorem{lem}[thm]{Lemma}
\newtheorem{prop}[thm]{Proposition}
\newtheorem{rem}[thm]{Remark}

\newcommand\JaL{Ja\-cquet-Lang\-lands}
\newcommand\JLc{\JaL\ cor\-res\-pon\-dence}

\def\subrep{sub-re\-pre\-sen\-ta\-tion}
\def\rep{re\-pre\-sen\-ta\-tion}

\newcommand\second{{$2$nd}}

\def\ie{{i.e.}}
\def\eg{{e.g.}}
\def\cf{{\it cf. \/}}

\def\lam{{\lambda}}
\def\uni{{\varpi}}

\def\color{\varrho}

\def\s{\sigma}

\def\Lam{{\Lambda}}

\def\A{{\mathbb A}}
\def\C{{\mathbb C}}
\def\F{{\mathbb F}}

\def\N{{\mathbb N}}
\def\Q{{\mathbb {Q}}}
\def\R{{\mathbb {R}}}

\def\Z{{\mathbb Z}}

\def\B{{\mathcal B}}

\def\GG{{\mathcal G}}
\def\I{{\mathcal I}}
\def\O{{\mathcal{O}}}

\def\Vals{{\mathcal{V}}}

\def\GammaU{{\Gamma_{\!\!\:U}}}

\def\opi{{\overline{\pi}}}

\DeclareMathAlphabet{\mathscr}{OT1}{pzc}{m}{it}

\newcommand\tensor[1][{}]{{\otimes_{#1}}}

\def\ra{{\rightarrow}}

\def\hra{{\,\hookrightarrow\,}}
\def\minusset{{-}}

\def\sub{\subseteq}
\def\sup{\supseteq}
\def\({\left(}
\def\){\right)}
\def\isom{{\;\cong\;}}
\def\normal{{\unlhd}}

\def\normali{{\lhd}}
\def\co{{\,{:}\,}}
\def\divides{{\,|\,}}
\newcommand\suchthat{{\,:\ \,}}
\newcommand\subjectto{{\,|\ }}

\newcommand\comp[1]{{{#1}^{\operatorname{c}}}}

\DeclareMathOperator{\spec}{spec}

\DeclareMathOperator{\Ker}{Ker}

\DeclareMathOperator{\diag}{diag}

\DeclareMathOperator{\Ind}{Ind}
\DeclareMathOperator{\Aut}{Aut}

\newcommand\JL{{\operatorname{JL}}}
\newcommand{\Sp}{\operatorname{Sp}}%
\DeclareMathOperator{\Cp}{C}

\DeclareMathOperator{\End}{End}

\DeclareMathOperator{\mychar}{char} \DeclareMathOperator{\Cent}{Z}

\newcommand{\Norm}[1][]{{\operatorname{N}_{#1}}}
\DeclareMathOperator{\Gal}{Gal}

\newcommand\cond[2][!]{{\operatorname{cond}_{\if!#1\relax\else{\comp{#1}}\fi}(#2)}}

\renewcommand\L[2]{{\operatorname{L}^{#1}(#2)}}
\newcommand\M[1][d]{{\operatorname{M}_{#1}}}
\newcommand\GL[1][d]{{\operatorname{GL}_{#1}}}
\newcommand\PGL[1][d]{{\operatorname{PGL}_{#1}}}
\newcommand\PGU[1][d]{{\operatorname{PU}_{#1}}}

\newcommand\SL[1][d]{{\operatorname{SL}_{#1}}}

\newcommand\PSL[1][d]{{\operatorname{PSL}_{#1}}}

\newcommand\abs[2][F]{|{#2}|_{#1}}
\newcommand{\norm}[1]{{\vert\vert {#1} \vert \vert}}

\newcommand{\set}[1]{{\{#1\}}}
\newcommand{\card}[1]{{\left|{#1}\right|}}
\newcommand\ideal[1]{{\left<{#1}\right>}}
\newcommand\sg[1]{{\ideal{#1}}}

\newcommand\dimcol[2]{{[{#1}\!:\!{#2}]}}

\newcommand\db[1]{{(\:\!\!({#1})\:\!\!)}}
\newcommand\mul[1]{{#1^{\times}}}
\newcommand\md[2][d]{{\mul{#2}/{\mul{#2}}^{#1}}}
\newcommand\restrict[2]{{{#1}|_{#2}}}
\newcommand\innprod[3][]{{\left<{#2},{#3}\right>_{#1}}}

\newcommand\He[2][G]{{\operatorname{H}(#1,#2)}}

\newcommand\valF{{\nu_0}}

\newcommand\Tref[1]{{Theorem \ref{#1}}}
\newcommand\Eref[1]{{Example \ref{#1}}}
\newcommand\Pref[1]{{Proposition \ref{#1}}}

\newcommand\Cref[1]{{Corollary \ref{#1}}}

\newcommand\Rref[1]{{Remark \ref{#1}}}

\newcommand\eq[1]{{(\ref{#1})}}

\newcommand\Eq[1]{{Equation \eq{#1}}}

\newcommand\defin[1]{{\it{#1}}}
\long\def\half#1\halved{{\footnotesize{#1}}}

\newcommand\dc[3]{{{#1}\backslash{#2}/{#3}}}
\newcommand\dom[2]{{{#1}\backslash{#2}}}

\newcommand\paper[6]{{{#1},\ {\it{#2}},\ {#3}\ {#4},\ {#5},\ ({#6}).}}
\newcommand\book[4]{{{#1},\ {{#2}},\ {#3},\ {#4}.}}
\newcommand\thesis[4]{{{#1},\ {{#2}}, Doctoral Dissertation, {#3},\ {#4}.}}

\newcommand\paperTrans[9]{{{#1},\ {\it{#2}},\ {#3}, {#4}, {#5}. Trans. {#6}, {#7}, {#8}, ({#9}).}}

\newcommand\absdot[1][]{\abs[#1]{\;\!\cdot\:\!}}
\newcommand\dd[4][!]{\abs[#2]{\:\!\det{\if!#1\relax\:\!\!\else(#1)\fi}}^{#3} #4}
\newcommand\HTplus{{\oplus}}

\newcommand\yy[1]{{{#1}_{{\:\!} 1 {\:\!\!} / {\;\!\!} y}}}
\def\ii{{\vec{\imath}\,}}
\def\jj{{\vec{\jmath}\,}}
\newcommand\automrep[1][G']{{\L2{\dom{{#1}(k)}{{#1}(\A)}}}}

\newcommand\binomq[3]{{\genfrac{[}{]}{0pt}{2}{#1}{#2}_{#3}}}
\newcommand\con[3]{{{#1}({#2},{#3})}}
\def\Aspec{{\mathfrak S}}

\def\id{{\operatorname{id}}}

\newcommand\Section[1]{\setcounter{equation}{0}\section{#1}}

\newif\ifXY
\XYfalse
\ifXY
\usepackage{xy}
\xyoption{all}
\fi

\begin{document}

\title
{Ramanujan Complexes of Type $\tilde{A_d}$}

\def\HUJI{Inst. of Math., Hebrew Univ., Givat Ram, Jerusalem 91904,
Israel}
\def\BIU{Dept. of Math., Bar-Ilan University, Ramat-Gan 52900,
Israel}
\def\YALE{Dept. of Math., Yale University, 10 Hillhouse Av., New-Haven CT
06520, USA}

\author{Alexander Lubotzky}
\address{\HUJI}
\email{alexlub@math.huji.ac.il}

\author{Beth Samuels}
\address{\YALE
}
\email{beth.samuels@yale.edu}

\author{Uzi Vishne }
\address{\YALE} \curraddr{\BIU}

\email{vishne@math.biu.ac.il}

\thanks{Partially supported by grants from NSF and BSF (U.S.-Israel).}

\renewcommand{\subjclassname}{
      \textup{2000} Mathematics Subject Classification}

\date{Submitted Sep. 1, 2003; revised Feb. 19, 2004}

\dedicatory{With love and admiration to Hillel F\"urstenberg, a
teacher and friend}

\begin{abstract}
We define and construct Ramanujan complexes. These are simplicial
complexes which are higher dimensional analogues of Ramanujan
graphs (constructed in \cite{LPS}). They are obtained as quotients
of the buildings of type $\tilde{A}_{d-1}$ associated with
$\PGL[d](F)$ where $F$ is a local field of positive
characteristic.
\end{abstract}

\maketitle

\section{Introduction}

A finite $k$-regular graph $X$ is called a Ramanujan graph if for
every eigenvalue $\lam$ of the adjacency matrix $A = A_X$ of $X$
either $\lam = \pm k$ or $\abs[]{\lam} \leq 2\sqrt{k-1}$. This
term was defined in \cite{LPS} where some explicit constructions
of such graphs were presented, see also \cite{Ma1},
\cite{alexbook}, \cite{Morg}.

These graphs were obtained as quotients of the $k$-regular tree $T
= T_k$, for $k = q+1$, $q$ a prime power, divided by the action of
congruence subgroups of $G = \PGL[2](F)$. Here $F$ is a
non-archimedean local field with residue field of order $q$, and
the tree $T$ is the Bruhat-Tits building associated with $G$,
which is a building of type $\tilde{A}_1$. The (proved) Ramanujan
conjecture for $\GL[2]$ was an essential ingredient in the proof
that the graphs are indeed Ramanujan, see \cite{alexbook}.

The number $2\sqrt{k-1}$ plays a special role in the definition of
Ramanujan graphs because of the Alon-Boppana theorem (see
\cite{LPS}), which proves that this is the best possible bound for
an infinite family of $k$-regular graphs. A conceptual explanation
was given by Greenberg \cite{Gr}, \cite[Thm.~4.2.7]{alexbook} (see
also
\cite{GZ}): for a connected graph $X$, let
$\rho(X)$ denote the norm of the adjacency operator $A$ on
$\L2{X}$ (so $\rho(T_k) = 2\sqrt{k-1}$\,); then, Greenberg showed
that no upper bound on the non-trivial eigenvalues of finite
quotients of $X$ is better then $\rho(X)$.

These considerations motivated Cartwright, Sol\'{e} and \.{Z}uk
\cite{CSZ} to suggest a generalization of the notion of Ramanujan
graphs from finite quotients of $T_k$ --- which is an
$\tilde{A}_1$ building --- to the simplicial complexes obtained as
finite quotients of $\B = \B_d(F)$, the Bruhat-Tits building of
type $\tilde{A}_{d-1}$ associated with the group $G = \PGL[d](F)$.
The vertices $\B^0$ of the building are labelled by a `color'
function $\color \co \B^0 \ra \Z/d\Z$, and we may look at the
$d-1$ colored adjacency operators $A_k$, $k = 1,\dots,d-1$ on
$\L2{\B^0}$, called the Hecke operators. They are defined by
\begin{equation}\label{Akdef}
(A_kf)(x) = \sum{f(y)}
\end{equation}
where the summation is over the neighbors $y$ of $x$ such that
$\color(y)-\color(x) = k$ in $\Z/d\Z$.

These operators $A_k$ are bounded, normal, and commute with each
other. Thus, they have a simultaneous spectral decomposition, and
the spectrum $\Aspec_d$ of $(A_1,\dots,A_{d-1})$ on $\L2{\B^0}$
was computed explicitly as a subset of $\C^{d-1}$ (see Subsection
\ref{ss:spec} below). This set is, of course, contained in the
Cartesian product $\Aspec_{d,1}\times \dots \times
\Aspec_{d,d-1}$, where $\Aspec_{d,k}$ is the spectrum of $A_k$,
but it is not equal to the product.

\begin{defn} [{following \cite{CSZ}}]
A finite quotient $X$ of $\B$ is called a Ramanujan complex if the
eigenvalues of every non-trivial simultanenous eigenvector $v \in
\L2{X}$, $A_k v = \lam_k v$, satisfy $(\lam_1,\dots,\lam_{d-1})
\in \Aspec_d$.
\end{defn}
(See Subsection \ref{ss:spec} for more detailed explanations, and
in particular for a description of the trivial eigenvalues. See
also \cite{JL} for a definition and construction of Ramanujan
complexes which are not simplicial).

Cartwright et al. \cite{CSZ} also suggested a way of obtaining such
Ramanujan complexes: assume $F$ is a local field of positive
characteristic; let $\Gamma$ be a cocompact arithmetic lattice of
$G = \PGL[d](F)$ of inner type, and $\Gamma(I)$ a congruence
subgroup of $\Gamma$. They conjectured that the quotients
$\dom{\Gamma(I)}{\B}$ are Ramanujan complexes.
The work of Lafforgue in the last few years, which proved the
Ramanujan conjecture for $\GL[d]$ in characteristic $p$ (an
extension of Drienfeld's work for $\GL[2]$ in characteristic $p$
and of Deligne's for $\GL[2]$ in characteristic zero) provided hope
that these combinatorial applications could be deduced.

The current work, which started from the challenge to prove the
conjecture in \cite{CSZ}, shows that for general $d$, the story is
more subtle. It turns out that most of these quotients are indeed Ramanujan,
but not all. To describe our results, let us first
introduce some notation.

Let $k$ be a global field of characteristic $p > 0$, and $D$, a
division algebra of degree $d$ over $k$. Denote by $G'$ the
$k$-algebraic group $\mul{D}/\mul{k}$, and fix a suitable
embedding of $G'$ as a linear group (see Section~\ref{inner}). Let
$T$ be the finite set of valuations of $k$ for which $D$ does not
split. We assume that for every $\nu \in T$, $D_{\nu} = D
\tensor[k] k_{\nu}$ is a
division algebra.
Let $\valF$ be a valuation of $k$ which is not
in $T$, and $F = k_{\valF}$. Let
\begin{equation}\label{R0def}
R_0 = \set{x \in k \suchthat
\nu(x)\geq 0 \quad \mbox{for every $\nu \neq \valF$}}.
\end{equation}
Then $\Gamma= G'(R_0)$ is a discrete subgroup of $G'(F)$, and the
latter is isomorphic to $G(F) = \PGL[d](F)$, as $F$ splits $D$. By
general results, $\Gamma$ is in fact a cocompact lattice in $G(F)$
--- an ``arithmetic lattice of inner type''. Let $\B = \B_d(F)$ be
the Bruhat-Tits building of $G(F)$, then $\B^0 \isom G(F)/K$,
where $K = G(\O)$ is a maximal compact subgroup ($\O$ is the ring
of integers in $F$). $G(F)$ acts on $\B$ by left translation.

For $0 \neq I \normali R_0$ an ideal (note that
$R_0$ is a principal ideal domain),
we have the principal congruence subgroup
\begin{equation}\label{condef}
\Gamma(I)  = \con{G'}{R_0}{I} = \Ker(G'(R_0) \ra G'(R_0/I)).
\end{equation}

{\emph{In the following two theorems we assume the global \JLc\
for function fields, see \Rref{valid} below regarding this
assumption.}}

\begin{thm}\label{IntroT1}
If $d$ is prime, then for every
$0 \neq I \normali R_0$, $\dom{\Gamma(I)}{\B}$ is a Ramanujan complex.
\end{thm}

So for $d$ prime, the Cartwright-Sol\'{e}-\.{Z}uk conjecture is indeed
true.
On the other hand, for general $d$:
\begin{thm}\label{IntroT2}
(a) For every $d$, if $I$ is prime to some valuation $\theta \in T$,
\ie\ $\theta(a) = 0$ for some $\theta \in T$ and some $a \in I$,
then $\dom{\Gamma(I)}{\B}$ is a Ramanujan complex.

(b) If $d$ is not a prime, then there exist (infinitely many) ideals $I$
such that $\dom{\Gamma(I)}{\B}$ are not Ramanujan.
\end{thm}

Theorem \ref{IntroT1} may suggest
that in positive characteristic, if $d$ is a prime,
then every finite quotient of $\B$ is
Ramanujan.
We do not know if this is indeed the case (which would
be truly remarkable), but at least in the zero characteristic analog
there are counter examples.
Indeed, in Section \ref{outer} we show that if $E$ is a
non-archimedean local field of characteristic zero,
then congruence quotients of $\B = \B_d(E)$ can be
non-Ramanujan for every $d \geq 4$. This happens if $\Gamma$ is
taken to be an arithmetic group of outer type.

\begin{thm}\label{IntroT3}
Let $E$ be a non-archimedean local field of characteristic zero,
and assume $d \geq 4$. Then $\B_d(E)$ has
infinitely many non-Ramanujan quotients.
\end{thm}
For a discussion of the case $d = 3$, see \cite{Ballantine}.
The proof of Theorem \ref{IntroT1} and \ref{IntroT2}(a) follows in
principle the line of proof for Ramanujan graphs, as in
\cite{alexbook}.
The problem is transferred to representation theory.

\begin{prop}\label{S=sA}
Let $\Gamma$ be a cocompact lattice in $G(F) = \PGL[d](F)$.
Then $\dom{\Gamma}{\B}$ is a Ramanujan complex iff
every
irreducible spherical infinite-dimensional \subrep{} of
$\L2{\dom{\Gamma}{G(F)}}$ is tempered.
\end{prop}

The strategy now is to start with an irreducible \subrep{} $\rho$
of $\L2{\dom{\Gamma(I)}{G'(F)}}$. By Strong Approximation, one can
show (see Subsection \ref{ss:SA} below) that $\rho$ is a local
factor of an ad\`{e}lic automorphic \rep{} $\pi' =
\tensor{\pi'_{\nu}}$ in $\L2{\dom{G'(k)}{G'(\A)}}$ such that
$\pi'_{\valF} = \rho$, where $\A$ is the ring of ad\`{e}les of
$k$.

We can view $\pi'$ as an automorphic
\rep{} of $\mul{D}(\A)$.
Then, the \JLc\ associates with $\pi'$ an automorphic \rep{}
$\pi = \tensor{\pi_{\nu}}$ in $\L2{\dom{\GL[d](k)}{\GL[d](\A)}}$,
such that $\pi_{\valF} = \pi'_{\valF}$. We then appeal to the work of
Lafforgue, who proved that if $\pi$ is cuspidal, then $\pi_{\nu}$
is tempered for every unramified
$\nu$, and in particular $\pi_{\valF} =
\rho$ is tempered.

Now, the cuspidality issue is exactly what distinguishes between
the cases where $d$ is a prime and where $d$ is a composite number.
If $d$ is prime, then all infinite-dimensional irreducible
\subrep{s} of $\L2{\dom{\GL[d](k)}{\GL[d](\A)}}$ are cuspidal (and
the others are one-dimensional, and are responsible for the
``trivial'' eigenvalues, see Subsection \ref{ss:spec}). Thus
Theorem \ref{IntroT1} can be proved.

On the other hand, when $d$ is not a prime, there is a ``residual
spectrum'' and $\pi$ may
be there, in which case $\pi_{\valF}$ is not tempered. Theorem
\ref{IntroT2} (both parts (a) and (b)) is proved by a careful
analysis of the image of the \JaL\ map, as described in \cite{HT}.

The proof of Theorem \ref{IntroT3} is different. We apply the
method of Burger-Li-Sarnak \cite{BLS1},\cite{BLS2} who showed how
the existence of large ``extended arithmetic subgroups'' in
$\Gamma(I)$ can affect the spectrum. For arithmetic lattices
associated to Hermitian forms (unlike the case of inner type),
such ``large'' subgroups do exist, but anisotropic Hermitian forms
(with enough variables) exist only if $\mychar(F) = 0$.

\begin{rem}\label{valid}
{\textnormal {The global \JLc\ is proved in the literature for
fields of characteristic zero (see \Tref{HTJL} below and
\cite[Thm.~VI.1.1]{HT}). It is likely that the theorem is valid in
exactly the same formulation in positive characteristic, and it
seems (to some experts we consulted) that a proof can be worked
out using existing knowledge. So far, this task has not been
carried out. We hope that our work will give some additional
motivation to complete this gap in the literature.
}}

{\textnormal {W.~Li \cite{Li2} managed to prove the existence of
Ramanujan complexes of type $\tilde{A}_{d}$ in positive
characteristic, avoiding the use of the \JLc, and in fact also not
using Lafforgue's theorem, appealing to \cite{LRS} instead.
In order to apply this method,
one needs the division algebra to be ramified in at least four
places, and therefore it does not cover the case of algebras ramified
in two places. This case is crucial for our next work, \cite{paperII},
in which we give an explicit construction of Ramanujan complexes.
On the other hand, we have to assume that in the ramification
points the algebra is completely ramified, while Li requires this
assumption in only two prime places.}}

{\textnormal {We recently learned that Alireza Sarveniazi
\cite{Ali} has also given a construction of Ramanujan complexes.}}
\end{rem}

The paper is organized as follows: in Section \ref{Prelim1} we
describe briefly the building $\B$, the operators $A_k$,
the local representation theory, and, in particular, we prove
\Pref{S=sA} above. In Section \ref{Prelim2} we show how strong
approximation enables one to pass from the local theory to the
global one.
In Section
\ref{Prelim3} we survey the global theory: Lafforgue's theorem,
the residual spectrum, and the \JLc. After the preparations we
prove Theorems \ref{IntroT1} and \ref{IntroT2} in Section
\ref{inner}, and Theorem \ref{IntroT3},
in Section \ref{outer}.

Much of the material of Sections \ref{Prelim1}--\ref{Prelim3} is
well known to experts, but since we expect (and hope) the paper
will have readers outside representation theory and automorphic
forms, we tried to present the material in a suitable way for
non-experts.

{\bf Acknowledgements}. We are indebted to
M.~Harris,
R.~Howe,
D.~Kazhdan,
E.~Lapid,
S.~Miller,
E.~Sayag,
I.~Piatetski-Shapiro,
T.~Steger,
M.F.~Vigneras,
and especially to
J.~Rogawski and J.~Cogdell
for many helpful discussions while working on the project. This work
was done while {the first named author visited and the third named
author held a post-doc position} at Yale, whose hospitality and
support are gratefully acknowledged. We also thank the NSF and the
BSF US-Israel for their support.

\Section{Affine buildings and \rep{s} of the local group}\label{Prelim1}

In this section, $F$ is a non-Archimedean local field of
arbitrary characteristic,
$\O$ its ring of integers, and $\uni \in \O$ a
uniformizer. Let $\valF \co F \ra \Z$ denote the valuation of $F$.

\subsection{Affine buildings of type $\tilde{A}_{d-1}$}

Recall that a complex is a structure composed of $i$-cells, where
the $0$-cells are called vertices, and every $i$-cell is a set of
$i+1$ vertices.
A complex is simplicial if every subset of a cell is also a cell.

We will now describe the affine building $\B = \B_d(F)$ associated to
$\PGL(F)$, which is an (infinite) simplicial complex.
Consider the $\O$-lattices of full rank in $F^d$.
We define an equivalence relation on lattices by
setting $L \sim s L$ for every $s \in \mul{F}$. Since $\mul{F}/\mul{\O}$ is
the infinite cyclic group generated by $\uni$, an equivalent definition
is that $L \sim \uni^{i} L$ for every $i \in \Z$.

By $\B^i$ we denote the set of $i$-cells of $\B$. The vertices
$\B^0$ are the equivalence classes of lattices. There is an edge
($1$-cell)$(x,x')$, from $x=[L]$ to $x'=[L'] \in \B^0$, if  $\uni
L \sub L' \sub L$. Notice that this is a symmetric relation, since
then $\uni L' \sub \uni L \sub L'$. The quotient $L/\uni L$ is a
vector space of dimension $d$ over the field $\O / \uni \O \isom
\F_{q}$.

As $i$-cells of $\B$ we take the complete subgraphs of size $i+1$
of $\B^0$. It immediately follows that $\B$ has $(d-1)$-cells
(corresponding to maximal flags in quotients $L/\uni L$). It also
follows that there are no higher dimensional cells. We call $L_0 =
\O^d \sub F^d$ the \defin{standard lattice}.

For every lattice $L$, there is some $i$ such that $\uni^i L \sub
L_0$ (it then follows that every two lattices of maximal rank are
commensurable). We define a color function $\color \co \B^0 \ra
\Z/d$, by
\begin{equation} \label{colordef}
\color(L) = \log_q \dimcol{L_0}{\uni^iL}
\end{equation}
for $i$ large enough; the color is well defined since
$\log_q\dimcol{\uni^iL}{\uni^{i+1}L} = d$. In a similar way, the
color of an ordered edge $(x,y) \in \B^1$ is defined to be
$\color(x)-\color(y) \pmod{d}$.

The group $\GL(F)$ acts on lattices by its action on bases; the
scalar matrices carry a lattice to an equivalent lattice, so $G =
\PGL(F)$ acts (transitively) on the vertices of $\B$. Since the
action of $\GL(F)$ preserves inclusion of lattices, $G$ respects
the structure of $\B$, and in particular the color of edges. Note
that $\GL(F)$ does not preserve the color of vertices, but
$\SL(F)$ does.

The stabilizer of $[L_0]$ is the maximal compact subgroup $K =
\PGL(\O)$. We can thus identify $\B^0$ with $G/K$, where $G$ acts
by multiplication from the left. A coset $gK \in G/K$ corresponds
to the lattice generated by the columns of
$g$ (so $[L_0]$ corresponds to the identity matrix).
The color of $gK$ can then be computed from the determinant of
$g$:
$$\det(g) \equiv \uni^{\color(gK)} \pmod{\mul{F}^d},$$
where $\mul{F}^d$ is the subgroup of $d$-powers in $\mul{F}$.

Let $\omega_k = \diag(\uni,\dots,\uni,1,\dots,1)$, where
$\det(\omega_k) = \uni^k$. The lattice corresponding to $\omega_k
K$ is obviously a neighbor of color $k$ of $[L_0]$. Let $\Omega_k$
be the set of neighbors of color $k$ of $[L_0]$.  Then $K$ acts
(as a subgroup of $G$) transitively on $\Omega_k$, so that $K
\omega_k K = \cup y K$, where the union is over $y K \in
\Omega_k$. Multiplying from the left by an arbitrary $g \in G$, we
see that the neighbors forming an edge of color $k$ with $gK$
are $\set{gyK}_{y K \in \Omega_k}$.

It follows that the operators $A_k$ (defined in \Eq{Akdef}) act
on functions $f \co G/K \ra \C$ by
$$(A_k f)(gK) = \sum_{y K \in \Omega_k} f(gyK) = \sum_{y K \in \Omega_k} \int_{y K} f(gx)dx = \int_{K \omega_k K} {f(gx) dx},$$
the integrals are normalized so that $\int_K dx = 1$.
See \cite{Mac} and \cite{Ballantine} for details.

\subsection{Spherical \rep{s}}\label{ss:sph}

In this section let $K = \GL(\O)$, which is a maximal compact subgroup of
$G = \GL(F)$.

As in \cite{alexbook},
we study the spectrum of the operators $A_k$ via representations
of $\GL(F)$. An irreducible admissible
representation of $G$ is called
\defin{$H$-spherical} if the \rep{} space has an $H$-fixed
vector, where $H \leq G$ is a subgroup. The $K$-spherical \rep{s}
are simply called \defin{spherical}. (A representation is
\defin{smooth} if every $v \in V$ is fixed under some open compact
subgroup, and \defin{admissible} if, moreover, the spaces fixed by
each open compact subgroup are finite dimensional).

The Hecke operators $A_k$ of the preceding subsection (defined in the
same way, as functions of $G/K$ for $G = \GL[d](F)$
rather than $G =\PGL[d](F)$) generate the Hecke
algebra $\He{K}$ of all bi-$K$-invariant
compactly supported
functions on $G$, with multiplication defined by
$$(A*A')(g) = \int_G A(x)A'(x^{-1}g)dx.$$
The $A_k$ commute with each other, and freely generate $\He{K}$ (\cf\ \cite[Sec.~V]{Mac}).

Let $\rho \co G \ra \End(V)$ be an admissible \rep{}; the Hecke
algebra acts on the representation space (see \cite[Eq.~(9)]{Cartier}) by
\begin{equation}\label{HGKact}
A\cdot v = \int_G A(x) (\rho(x))(v)\,dx,
\end{equation}
which is an integration over a compact set since $A$ is compactly
supported. It projects $V$ to the $K$-fixed subspace $V^K$ (which
is finite dimensional as the \rep{} is admissible). Moreover, if
$V$ is an irreducible $G$-module, then $V^K$ is an irreducible
$\He{K}$-module. Since $\He{K}$ is commutative and finitely
generated,
$V^K$ is one-dimensional in this case,
and consequently, every $v \in V^K$ is an eigenvector of all the $A_k$.

We describe how spherical representations are parameterized by
$d$-tuples of complex numbers, called the Satake parameters. For
details, the reader is referred to \cite{Cartier}.

Let $B$ denote a Borel subgroup of $G$ (\eg\ the upper triangular
matrices), $U$ its unipotent radical, and $T \isom B/U \isom
(\mul{F})^d$ a maximal torus of $B$. We then have $B = UT$ and $G
= BK = UTK$.

Recall that $\mul{F}/\mul{\O} = \sg{\uni}$. A
character $\chi \co \mul{F} \ra \mul{\C}$ is spherical if it is
trivial on the maximal compact subgroup of $\mul{F}$, namely
$\mul{\O}$. Such a character is, thus, determined by $z =
\chi(\uni)$, which is an arbitrary complex number. The character
is called \defin{unitary} iff $z \in S^1 = \set{w \in \C \suchthat \abs[]{w} = 1}$.

Every character $\chi \co T \ra \mul{\C}$ can be written as
$\chi(\diag(t_1,\dots,t_d)) = \chi_1(t_1)\dots\chi_d(t_d)$, for
characters $\chi_i \co \mul{F} \ra \mul{\C}$. $\chi$ is said to be
\defin{unramified} if the $\chi_i$ are spherical. Since $T \isom
B/U$, $\chi$ extends to a character of $B$. The symmetric group
$S_{d}$ acts on the characters by permuting the $\chi_i$.

The unitary induction of \rep{s}
from $B$ to $G$ is defined using the \defin{modular function}
\begin{equation}\label{modularP}
\Delta(b) = \abs{a_1}^{d-1}\abs{a_2}^{d-3}\dots
\abs{a_d}^{1-d},\qquad b\in B \end{equation} where
$a_1,a_2,\dots,a_d$ are the entries on the diagonal of $b$ (in
that order), and $\absdot[F]$ is the absolute value function of
$F$, normalized so that $\abs{x} = q^{-\nu_0(x)}$ where $q =
\card{\O/\uni\O}$. The induced representation $I_{\chi} =
\Ind_B^G(\chi)$, is the space of locally constant functions $f \co
G\ra \C$ such that
$$f(bg) = \Delta^{1/2}(b)\chi(b) f(g) ,\qquad b\in B,\, g\in G$$
with the action of $G$ from the right (by $g \cdot f(x) = f(xg)$).
The inclusion of the modular function $\Delta$ guarantees that if
$\chi$ is unitary, then there is an inner product $\innprod{f}{f'}
= \int_K{f(x)\overline{f'(x)}dx}$ on $I_{\chi}$, for which the
action of $G$ is unitary (these are called the spherical principal
series \rep{s}). However, the space $I_{\chi}$ still can be
unitary even if $\chi$ is not unitary (these are called spherical
complementary series \rep{s}), see Subsection \ref{ss:bounds}.

We remark that $I_{\chi}$ need not be irreducible. Two spaces $I_{\chi}$ and
$I_{\chi'}$ are isomorphic iff $\chi' = w\chi$ for some $w \in S_d$
(\cite[Subsec.~3.3]{Cartier}, \cite[Sec.~2.6]{Bump}).

Notice that $\absdot$ is spherical, so the modular function
$\Delta$ is an unramified character. If $g \in B \cap K$ then $g$
is upper triangular, with its diagonal entries invertible in $\O$.
Since $G = BK$ and $\chi$ is unramified, it follows that
\begin{equation}\label{fchidef}
f_\chi(bk) := \Delta^{1/2}(bk) \chi(bk)= \Delta^{1/2}(b) \chi(b),
\qquad b \in B,\, k\in K,
\end{equation}
is a well defined $K$-fixed function (
unique in $I_\chi$), which makes the induced \rep{} $\rho \co G
\ra \End(I_{\chi})$ spherical. By definition, $\rho$ is determined
by the numbers $z_i = \chi_i(\uni) =
\chi(\diag(1,\dots,1,\uni,1,\dots,1))$, called the Satake
parameters of $\chi$, where $\chi_i$ are the diagonal components
of $\chi$, which is a \subrep{} of
$\restrict{\Delta^{-1/2}\rho}{B}$. The \rep{s} which are well
defined on $\PGL(F)$ are those with $z_1 \cdots z_d = 1$ (since
they need to be trivial on the center of $\GL[d](F)$).

Let $\s_{k}(z_1,\dots,z_d)$ be the $k$th elementary symmetric
function, \ie\ $\s_k(z_1,\dots,z_d) =
\sum_{i_{1}<\dots<i_{k}}z_{i_1}\dots z_{i_k}$.

\begin{prop}\label{fchi}
The function $f_{\chi}$ is an eigenfunction of the $A_k$,
$A_k f_{\chi} = \lam_k f_{\chi}$,
where $\lam_k = q^{k(d-k)/2} \s_k(z_1,\dots,z_d)$.
\end{prop}
\begin{proof}
Since $\He{K}$ acts on $I_{\chi}$ and preserves the $K$-fixed subspace
$\sg{f_{\chi}}$, $f_{\chi}$ is an eigenvector of the $A_k$.

It is enough to compute $A_kf_{\chi}$ at the point $g = 1$ (noting that
$f_{\chi}(1) = 1$). For every subset $C \sub \set{1,\dots,d}$ of size
$k$, let $\Omega_{k,C}$ be the set of upper triangular matrices $m$
such that $m_{ii} = \uni$ if $i \in C$, $m_{ii} = 1$ if $i \not
\in C$, $m_{ij}$ is in some fixed lifting of $\O/\uni\O$ to $\O$
if $i \in C$ and $j \not \in C$, and $m_{ij} = 0$ otherwise.
For example, for $d = 4$ and $k = 2$ the sets are
\newcommand\FMatrix[4]{{\(\!\!\! \begin{array}{cccc}#1 \\[-0.1cm] #2 \\[-0.1cm] #3 \\[-0.1cm] #4 \end{array}\!\!\! \)\! }}
\newcommand\FLine[4]{{#1}\!\! &\!\! {#2}\!\! &\!\! {#3}\!\! &\!\! {#4}}
{\small
$$
\FMatrix
{\FLine{\uni}{0}{*}{*}}
{\FLine{0}{\uni}{*}{*}}
{\FLine{0}{0}{1}{0}}
{\FLine{0}{0}{0}{1}},
\FMatrix
{\FLine{\uni}{*}{0}{*}}
{\FLine{0}{1}{0}{0}}
{\FLine{0}{0}{\uni}{*}}
{\FLine{0}{0}{0}{1}},
\FMatrix
{\FLine{\uni}{*}{*}{0}}
{\FLine{0}{1}{0}{0}}
{\FLine{0}{0}{1}{0}}
{\FLine{0}{0}{0}{\uni}},
\FMatrix
{\FLine{1}{0}{0}{0}}
{\FLine{0}{\uni}{0}{*}}
{\FLine{0}{0}{\uni}{*}}
{\FLine{0}{0}{0}{1}},
\FMatrix
{\FLine{1}{0}{0}{0}}
{\FLine{0}{\uni}{*}{0}}
{\FLine{0}{0}{1}{0}}
{\FLine{0}{0}{0}{\uni}},
\FMatrix
{\FLine{1}{0}{0}{0}}
{\FLine{0}{1}{0}{0}}
{\FLine{0}{0}{\uni}{0}}
{\FLine{0}{0}{0}{\uni}}.
$$
}

There is a one-to-one correspondence between neighbors $y K \in \Omega_k$
of $[L_0]$ and subspaces of co-dimension $k$ of $L_0/\uni L_0$,
so using the correspondence between
matrices and lattices mentioned above, $\Omega_k = \cup_C \Omega_{k,C}$.

Fix a subset $C$, and let $s = \sum_{i\in C} i$. The number of
matrices in $\Omega_{k,C}$ is $q^{(d-k+1)+\dots+d-s}$, while
$\Delta^{1/2}(y) = \Pi_{i\in C} \abs{\uni}^{(d+1)/2-i} =
q^{-k(d+1)/2+s}$ for every $yK \in \Omega_{k,C}$. It follows that
the sum of $\Delta^{1/2}(y)\chi(y)$ over $yK \in \Omega_{k,C}$ is
$q^{k(d-k)/2}\chi(y) = q^{k(d-k)/2}\Pi_{i\in C}z_{i}$, so by
summing
over all $C$ we obtain
$$(A_k f_{\chi})(1) =
\sum_{y K \in \Omega_k} \Delta^{1/2}(y)\chi(y)  = q^{k(d-k)/2}
\s_k(z_1,\dots,z_d).$$
\end{proof}

Now let $\rho \co G \ra \End(V)$ be an irreducible spherical
representation of $G=\GL(F)$ with a unique (up to scalar
multiples) $K$-invariant vector $v_0 \in V$.  Let $\hat{V}$ be the
\rep{} contragredient to $V$. The space $\hat{V}^K$ is dual to
$V^K$ and thus one-dimensional. Choose $\hat{v_0}\in \hat{V}^K$
such that $\innprod{v_0}{\hat{v_0}} = 1$ (where $\innprod{\,}{\,}$
is the action of $\hat{V}$ on $V$). Define the bi-$K$-invariant
function $\psi(g) = \innprod{\rho(g)v_0}{\hat{v_0}}$ (so that
$\psi(1) = 1$).  $\psi(g)$ is called a \defin{spherical function}.
As explained earlier, $v_0$ is an eigenvector of $\He{K}$.  The
action of $\He{K}$ on $V^K = \sg{v_0}$ defines a homomorphism
$\omega \co \He{K} \ra \C$ by
\begin{equation}\label{Aact}
A\cdot v_0 = \omega(A)v_0.
\end{equation}

The action of $G$ (from the right)on the space of functions
$\set{f \co  G\ra \C}$ induces an action of the Hecke algebra on
this space, and by \Eq{HGKact} and the definition of $\psi$, we
find that
$A \cdot \psi = \omega(A) \psi$.
Using \Eq{HGKact} one can check that
\begin{equation}\label{psiomega}
\psi(g) = \omega(1_{KgK}) / \mu(KgK),
\end{equation}
where $1_{KgK}$ is the characteristic function
of $KgK \sub G$, and $\mu$ is the normalized measure.

\begin{prop} \label{omegapsi}
Using the notation as above, if $\rho_1$ and $\rho_2$ are
irreducible spherical representations, then $\rho_1 \cong
\rho_2$ iff $\psi_1 = \psi_2$ iff $\omega_1 = \omega_2$.
\end{prop}
\begin{proof}
Let $v_{1_0}$ and $v_{2_0}$ be the (unique) $K$-fixed vectors of
$\rho_1$ and $\rho_2$. The equivalence of $\psi_1 = \psi_2$ and
$\omega_1 = \omega_2$ follows at once from \Eq{psiomega}. If
$\rho_1 \cong \rho_2$ then it is obvious that $\psi_1 =
\innprod{\rho_1(\cdot)v_{1_0}}{\hat{v_{1_0}}} =
\innprod{\rho_2(\cdot)v_{2_0}}{\hat{v_{2_0}}}
 = \psi_2$. In the other
direction, assume $\psi_1 = \psi_2$, and let $V_i$ be the representation space
of $\rho_i$. Define a map from $V_i$
to $\sg{\psi_i G}$, the representation spanned by $\psi_i$, by
sending $v \in V_i$ to the function $g \mapsto \innprod{\rho(g)
v}{\hat{v_{i_0}}}$, for $i=1,2$. This is easily seen to be a
non-zero homomorphism (since $v_{i_0} \mapsto \psi_i \neq 0$),
which is an isomorphism since $V_i$ is irreducible. Then
$\rho_1 \isom \sg{\psi_1 G} = \sg{\psi_2 G} \isom \rho_2$.
\end{proof}

As a $G$-module, $I_{\chi}$ has finite composition length, so
it has only finitely many irreducible subquotients.

\begin{prop}[\cite{Cartier}]\label{sphrep}
Every irreducible spherical \rep{}
of $\GL[d](F)$ is isomorphic to
a subquotient of $I_\chi$ for some unramified character
$\chi$, which is unique up to permutation.
\end{prop}
\begin{proof}
Let $\rho \co G \ra \End(V)$ be an irreducible spherical
representation of $G$ with $v_0$ as its $K$-invariant
vector. Let $\omega \co \He{K} \ra \C$ be its corresponding
homomorphism, defined by \Eq{Aact}.

It can be shown \cite[Cor.~4.2]{Cartier}
that every such homomorphism
is of the form
$$\omega_\chi(A) = \int_{G} {A(x) f_{\chi}(x)\,dx}$$
for some unramified character
$\chi \co T \ra \C$ (unique up to
permutation) where $f_{\chi}$ is defined in \Eq{fchidef}.
Then
$$(A\cdot f_{\chi})(1) = \int_{G}{A(x)f_{\chi}(x) dx} = \omega_{\chi}(A)
= \omega_\chi(A) \cdot f_{\chi}(1),$$ and since $f_{\chi}$ is an
eigenvector, $A\cdot f_{\chi} = \omega_{\chi}(A)f_{\chi}$ for
every $A \in \He{K}$. Let $W$ be an irreducible subquotient of
$I_{\chi}$ in which $f_{\chi}$ has a non-zero image. By the
previous proposition $\rho$ is isomorphic to $W$, since
$\omega_{\chi} = \omega$.
\end{proof}

Thus, every spherical representation is determined by the Satake
parameters $z_i = \chi_i(\uni) =
\chi(\diag(1,\dots,1,\uni,1,\dots,1))$, for some unramified
$\chi$, uniquely determined up to permutation.

\begin{prop}\label{trade}
Let $f  \co G/K \ra \C$
be a simultaneous eigenvector of $A_1,\dots,A_{d-1}$. Then there
is an unramified character $\chi$ such that $f_{\chi}$ has
the same eigenvalues.
\end{prop}
\begin{proof}
Consider $f \co G \ra \C$, which is invariant with respect to $K$.
Let $\sg{fG}$ denote the linear span of the $G$-orbits of $f$,
where $G$ acts from the right. Taking this space modulo a maximal
sub-module not containing $f$, we obtain an irreducible spherical
representation, where $f$ is a (unique) $K$-fixed vector. By the
previous proposition, it is isomorphic to a subquotient of
$I_{\chi}$ for an unramified character $\chi$, where $f_{\chi}$ is
the unique $K$-fixed vector. By \Pref{omegapsi}, since the two
representation spaces are isomorphic, they induce the same
homomorphism $\omega \co \He{K} \ra \C$, namely $A_k\cdot f =
\omega(A_k)f$ and $A_k\cdot f_{\chi} = \omega(A_k)f_{\chi}$.
\end{proof}

Let $\rho \co G \ra \End(V)$ be a unitary representation, and
$\innprod{\,}{\,}$ the inner product defined on $V$. The functions
of the form
$$\rho_{v,w} \co g \mapsto \innprod{\rho(g)v}{w}$$
where $v,w \in V$ are called the
\defin{matrix coefficients} of $\rho$. Notice that if $V$ has a
$K$-fixed vector $v_0$ and $\innprod{v_0}{v_0} = 1$, then
$\rho_{v_0,v_0}$ is a spherical function. In the special case of
$I_{\chi}$, $\rho_{f_{\chi},f_{\chi}}(g) =
\int_K{f_{\chi}(xg)\,dx}$.

If $V$ is irreducible then
fixing $w\neq 0$, the map $v \mapsto \rho_{v,w}$ is an isomorphism of \rep{s}
(where $G$ acts on the space of functions from the right).

A \rep{} is called \defin{tempered}, if for some $0 \neq v,w \in
V$, $\rho_{v,w} \in \L{2+\epsilon}{G}$ for every $\epsilon > 0$.
The following equivalence is well known:

\begin{prop}\label{t=S}
An irreducible spherical unitary \rep{} is tempered iff its Satake
parameters have absolute value $1$.
\end{prop}

\subsection{Ramanujan complexes and the spectrum of $A_k$}\label{ss:spec}

Let $\Gamma$ be a cocompact lattice of $G = \PGL[d](F)$. Then
$\Gamma$ acts on $\B = G/K$ by left translation, and
$\dom{\Gamma}{\B}$ is a finite complex. The color function defined
on $\B^0$ (\Eq{colordef}) may not be preserved by the map $\B \ra
\dom{\Gamma}{\B}$. However, the colors defined on $\B^1$ by
$\color(x,y) = \color(x)-\color(y) \pmod{d}$ are preserved, since
they are determined by the index of (a representative of) $y$ as a
sublattice in (a representative of) $x$.

Since the Hecke algebra $\He[G]{K}$ acts on $G$ from the right,
and $\Gamma$ is acting from the left, the operators $A_k$ on
$\L2{\B} = \L2{G/K}$ induce colored adjacency operators on
$\dom{\Gamma}{\B}$.

It should be noted that if $\Gamma$ is torsion free, then $\gamma
x \neq x$ for any $\gamma \neq 1$ and any $x \in \B^0$,
so the underlying graph of $\dom{\Gamma}{\B}$ is simple.
Every cocompact lattice
has a finite index torsion free subgroup.

The trivial eigenvectors appear in $\L2{\dom{\Gamma}{\B}}$ but not
in $\L2{\B}$, since the former complex is finite. The trivial
eigenvectors can be constructed as follows.
The trivial \rep{} of $G$ is obviously spherical. Taking $\chi =
\Delta^{-1/2}$, we see that $f_{\chi}(g) = 1$ for every $g$ (see
\Eq{fchidef}), and the action of $G$ on the subspace $\C f_{\chi}
\sub I_{\chi}$ is trivial. The Satake parameters of the trivial
representation are thus $z_i =
\Delta^{-1/2}(\diag(1,\dots,1,\uni,1,\dots,1)) = q^{-(d-2i+1)/2}$.
More generally, since $G(F)/\PSL[d](F)K \isom \md[d]{F}$, $G$ has
$d$ one-dimensional spherical representations. Fixing $\zeta$ such
that $\zeta^d=1$,
\begin{equation}\label{trivdef}
\chi(g) = \Delta^{-1/2}(g)\zeta^{\valF(\det(g))}
\end{equation}
corresponds to a one-dimensional \rep{}, with the $K$-fixed vector
$f_\chi(g) = \zeta^{\valF(\det(g))}$ (since $\det(k) \in \mul{\O}$
for every $k \in K$). The Satake parameters in this case are
$\zeta q^{-(d-1)/2},\dots,\zeta q^{(d-1)/2}$, and the eigenvalues
are
$$\zeta^k \cdot q^{k(d-k)/2}\s_k(q^{-(d-1)/2},q^{-(d-3)/2},\dots,q^{(d-1)/2}).$$
Let $t = \dimcol{\Gamma}{\Gamma \cap \PSL[d](F)}$. By
\Eq{trivdef}, the trivial eigenvector $f_{\chi}$ is well defined
on $\dom{\Gamma}{\B}$ iff $\zeta^{d/t} = 1$, and the respective
$d/t$ eigenvectors give rise to the  \defin{trivial} eigenvalues.
Since $G$ is infinite, these $f_{\chi}$ do not belong to
$\L2{\B}$.

\medskip

Let $\Aspec_{d,k} \sub \C$ denote the spectrum of the operator
$A_k$ acting on $\L2{\B}$.

\begin{defn}
The complex $\dom{\Gamma}{\B}$ is pseudo-Ramanujan if for each $k
= 1,\dots,d-1$, the non-trivial eigenvalues of $A_k$ acting on
$\L2{\dom{\Gamma}{\B}}$ belong to $\Aspec_{d,k}$.
\end{defn}

Let $\Aspec_d \sub \C^{d-1}$ denote the simultaneous spectrum of
$(A_1,\dots,A_{d-1})$ acting on the space $\L2{\B}$, namely, the
set of $(\lam_1,\dots,\lam_{d-1}) \in \C^{d-1}$ for which there
exist a sequence of unit vectors $v_n \in \L2{\B}$ such that
$\lim_{n \ra \infty} (A_k v_n - \lam_k v_n) = 0$ for every $k =
1,\dots, d-1$.

\begin{defn}
The complex $\dom{\Gamma}{\B}$ is Ramanujan if for every
non-trivial simultaneous eigenvector $f \in \L2{\dom{\Gamma}{\B}}$
of the $A_k$, the eigenvalues $(\lam_1,\dots,\lam_{d-1})$ belong
to $\Aspec_d$.
\end{defn}

Since the $A_k$ commute,\label{statement1} every eigenvalue of
$A_k$ can be obtained by a simultaneous eigenvector. Hence, a
Ramanujan complex is pseudo-Ramanujan. On the other hand,
$\Aspec_d$ is not the Cartesian product of the $\Aspec_{d,k}$. For
example, inverting the direction of edges in $\B$ carries $A_k$ to
$A_{d-k}$, so the operators $A_k$ and $A_{d-k}$ are adjoint to
each other. In particular for every $(\lam_1,\dots,\lam_d) \in
\Aspec_d$ we have that $\lam_{d-k} = \bar{\lam_k}$.

\begin{rem}\label{CanS}
The spectrum $\Aspec_{d,k}$ of $A_k$  on $\L2{\B}$ is equal to the
projection of $\Aspec_d$ on the $k$th component.
\end{rem}

\begin{rem}
If $d = 2$ or $d = 3$, then $\dom{\Gamma}{\B}$ is Ramanujan iff it
is pseudo-Ramanujan (indeed, for $d = 2$ the definitions coincide,
and for $d = 3$, $A_2$ is the adjoint operator of $A_1$).
\end{rem}

Let $S = \set{(z_1,\dots,z_d) \suchthat \abs[]{z_i} =1,\, z_1 \cdots z_d = 1}$
and $\sigma \co S \ra \C^{d-1}$ be the map defined by
$(z_1,\dots,z_d) \mapsto (\lam_1,\dots,\lam_{d-1})$, where
$$\lam_k = q^{k(d-k)/2}\s_k(z_1,\dots,z_d).$$

The theorem below is proved in \cite{Cw}.
For completeness, we sketch the proof,
following ideas from \cite{CM} (where the result was proved for $d = 3$).
First, we will need an easy calculus lemma:
\begin{lem}
Let $(a_n),(b_n)$ be positive series. If $\limsup(a_n
b_n^{2+\epsilon})\leq 1$ for every $\epsilon>0$ and $\set{a_n}$ is
bounded, then $\limsup(a_n b_n^{2})\leq 1$.
\end{lem}
\begin{proof}
Otherwise let $C > 1$ be an upper bound of $\set{a_n}$ and $p =
\limsup(a_n b_n^2)>1$, and take $\epsilon < 2 \log(p)/\log(C)$.
Then $$p = \limsup (a_n^{\frac{\epsilon}{2+\epsilon}}
a_n^{\frac{2}{2+\epsilon}} b_n^2) \leq
C^{\frac{\epsilon}{2+\epsilon}} \limsup(a_n
b_n^{2+\epsilon})^{\frac{2}{2+\epsilon}} < C^{\epsilon/2},$$ a
contradiction.
\end{proof}

\begin{thm}\label{specAk}
The simultaneous spectrum $\Aspec_d$ is equal to $\sigma(S)$.
\end{thm}
\begin{proof}
Let $\underline{z} = (z_1,\dots,z_d) \in S$. Then the
corresponding character $\chi$ is unitary, and the irreducible
subquotient generated by $f_{\chi}$ of the induced \rep{}
$I_{\chi}$, is tempered (\Pref{t=S}). Thus, the corresponding
spherical function $\psi_{\chi}$ is in $\L{2+\epsilon}{G}$ for
every $\epsilon > 0$. We already saw that $\psi_{\chi}$ is an
eigenvector, however it does not belong to $\L2{G}$. In order to
show that $\sigma(\underline{z})$ is in the spectrum, we twist
$\psi_{\chi}$ to elements of $\L2{G}$ which are "almost"
eigenvectors, and their almost-eigenvalues converge to
$\sigma(\underline{z})$.

By \Pref{fchi} (and since $\psi_{\chi}$ is the spherical function
associated to $f_{\chi}$), $A_k \psi_{\chi} = \lam_k \psi_{\chi}$
where $\lam_k = q^{k(d-k)/2}\s_k(z_1,\dots,z_d)$. For every vertex
$x \in \B^0 = G/K$, let $w(x)$ denote the distance (in $\B^1$) of
$x$ from the origin $[L_0]$.

Recall \cite[V.(2.2)]{Mac} that every double coset in
$\dc{K}{G}{K}$ has a unique representative of the form
$\diag(\uni^{\ell_1},\dots,\uni^{\ell_d})$ where $\ell_1 \geq
\dots \geq \ell_{d-1} \geq \ell_d = 0$; we call this
representative of $KgK$ the \defin{type} of $gK$, and note that
the number of vertices of this type is $\mu(KgK)$. Its distance
from $[L_0]$ is equal to $\ell_1$, so there are
$\binom{n+d-2}{d-2}< (n+d)^d$ types of distance $n$.

For $\delta > 0$, define a function $\psi_{\chi}^{\delta}$ on
$\B^0$ by $\psi_{\chi}^{\delta}(x) =
(1-\delta)^{w(x)}\psi_{\chi}(x)$. For each $n$, let $g_n$ denote
the type of the vertex of distance $n$ for which
$\mu(Kg_nK)\abs[]{\psi_{\chi}(g_nK)}^2$ is maximal.

To see that $\psi_{\chi}^{\delta} \in \L2{\B^0}$, compute that
\begin{eqnarray*} \sum_{x \in
\B^0}(1-\delta)^{2w(x)}\abs[]{\psi_{\chi}(x)}^2 & = &
\sum_{n=0}^{\infty}{(1-\delta)^{2n}\sum_{w(x) =
n}\abs[]{\psi_{\chi}(x)}^2} \\
& \leq & \sum_{n=0}^{\infty}{(1-\delta)^{2n}(n+d)^d
\mu(Kg_nK)\abs[]{\psi_{\chi}(g_nK)}^2}, \end{eqnarray*} and the
convergence follows from the root test once we show that
$\limsup{(\mu(Kg_nK)\abs[]{\psi_{\chi}(g_nK)}^2)^{1/n}}\leq 1$.
But since $\psi_{\chi} \in \L{2+\epsilon}{\B^0}$ for every
$\epsilon
> 0$, we have
$$\limsup{(\mu(Kg_nK)\abs[]{\psi_{\chi}(g_nK)}^{2+\epsilon})^{1/n}}\leq
1,$$ and the result follows from $\mu(Kg_nK)^{1/n} \leq q^{d}$ by
the lemma.

By the definition of $A_k$, $A_k \psi_{\chi}^{\delta}(x)$ is a sum
of $\psi_{\chi}^{\delta}(y)$ for neighbors $y$ of $x$, and the
distance of neighbors satisfies $\abs[]{w(y)-w(x)} \leq 1$. Since
$A_k \psi_{\chi} - \lam_k \psi_{\chi} = 0$, it follows that for
some constant $c$, $\norm{A_k \psi_{\chi}^{\delta} - \lam_k
\psi_{\chi}^{\delta}} \leq c \delta \norm{\psi_{\chi}^{\delta}}$
for every $\frac{1}{2} > \delta > 0$, showing that
$(\lam_1,\dots,\lam_{d-1})
\in \Aspec_d$.

Now let $(\lam_1,\dots,\lam_{d-1}) \in \Aspec_d
$, and let $z_1,\dots,z_d \in \C$ be numbers satisfying
$q^{k(d-k)/2}\s_k(z_1,\dots,z_d) = \lam_k$, with the added
property that $z_1 \dots z_d = 1$ (the $z_i$ are unique up to
order). We need to show that $(z_1,\dots,z_d) \in S$, implying
$(\lam_1,\dots,\lam_{d-1}) \in \sigma(S)$.

Let $v_n \in \L2{\B^0}$ be unit vectors such that $A_k v_n -
\lam_k v_n \ra 0$, for all $k$, and define a homomorphism $\omega
\co \He{K} \ra \C$ by $\norm{A v_n - \omega(A) v_n} \ra 0$ (here
we use the fact that the $A_k$ are bounded and generate $\He{K}$).
Then $\omega$ is continuous in the norm of the operators on
$\L2{\B^0}$, and in particular $\abs[]{\omega(A)} \leq \norm{A}$
for every $A \in \He{K}$ (otherwise take $\epsilon$ such that
$\abs[]{\omega(A)} > \norm{A}+\epsilon$, then
$(\frac{A}{\norm{A}+\epsilon})^n$ converges to zero but
$\omega((\frac{A}{\norm{A}+\epsilon})^n)$ does not).

For $\ell \geq 1$, let $H_\ell \in \He{K}$ be the characteristic
function of $K \diag(\uni^{d\ell},1,\dots,1) K$. We show that
while $\omega(H_\ell)$ is a certain combination of $z_r^{-d
\ell}$, the bound $\norm{H_{\ell}}$ is polynomial in $\ell$, thus
implying that $\abs[]{z_r} \geq 1$ for every $r$.

The vector $\psi_1$ associated to the trivial character $\chi = 1$
is strictly positive (since $f_1(x) > 0$ for every $x \in G/K$ and
$\psi_1(x) = \int_{K}{f_1(kx)dk}$), so if $H_\ell \psi_1 = b
\psi_1$ and $H_{\ell}^* \psi_1 = b' \psi_1$, we have
$\norm{H_{\ell}} \leq \sqrt{bb'}$ by Schur's criterion
\cite[p.~102]{Ped}. Let $p = ((d-1)\ell,-\ell,\dots,-\ell)$.

{}From
\cite[(3.5)]{Mac} and \cite[(3.3)]{Mac}, and using the limit
\begin{equation}\label{limit}
\lim_{(x_1,\dots,x_d) \ra
(1,\dots,1)}{\sum_{k=1}^d{\frac{x_k^m}{\prod_{i \neq k}{(x_k - x_i)}}}} = \binom{m}{d-1},
\end{equation}
we obtain $bb' =
(1-q^{-1})^{2(d-1)}\binom{d\ell}{d-1}\binom{d(\ell+1)-2}{d-1}q^{d(d-1)\ell}
< (d\ell)^{2d}q^{d(d-1)\ell} $, so $\norm{H_\ell} < (d\ell)^dq^{d(d-1)\ell/2}$.
(Note that the action of $\He{K}$ on the spherical functions in
\cite{Mac} is via the multiplication of the Hecke algebra, unlike
ours; see \Eq{HGKact}).

In a similar manner, $\omega(H_\ell)$ is equal to
$\widehat{c_{-p}}(\omega_s)$ of \cite[(3.3)]{Mac}, and has the
form $q^{d(d-1)\ell/2} \sum_{r = 1}^{d} {\alpha_r {z_r^{-
d\ell}}}$ where $\alpha_r = \prod_{i \neq r}\frac{z_i -
q^{-1}z_r}{z_i-z_r}$ if all the $z_i$ are different (see
\cite[III.(2.2)]{Mac}). {}From the continuity of $\omega$ we
proved
$$\sum_{r = 1}^{d} {\alpha_r {z_r^{- d\ell}}}  \leq C \ell^d$$
for some constant $C$ and every $\ell$.

Order the $z_i$ by absolute value, so that $\abs[]{z_1} \leq \dots
\leq \abs[]{z_d}$. Then $\alpha_1 \neq 0$ and from the last bound
it follows that $\abs[]{z_1}\geq 1$, but $z_1 \dots z_d = 1$ so
$(z_1,\dots,z_d) \in S$. If the $z_i$ are not assumed to be
different, one computes the coefficients of the $z_i^{-\ell d}$ by
\Eq{limit}, and the same arguments apply.
\end{proof}

The sets $\Aspec_{d,k}$ are explicitly described in \cite{CS2}:
$\Aspec_{d,k}$ is the simply connected domain with boundary the
complex curve
$$\set{q^{k(d-k)/2}\s_k(e^{i \theta},\dots,e^{i \theta}, e^{-(d-1)i\theta}) \suchthat \theta \in [0,2\pi]}$$
where $i = \sqrt{-1}$.

Notice that the equations $\lam_k =
q^{k(d-k)/2}\s_k(z_1,\dots,z_d)$ always have a solution, but
unless $(\lam_1,\dots,\lam_d) \in \Aspec_d$, the $z_i$ do not have
to be unitary---even if each $\lam_k \in \Aspec_{d,k}$.

In terms of characters, Theorem \ref{specAk} implies that the
eigenvalues corresponding to $f_{\chi}$ (see \Pref{fchi}) are in
the simultaneous spectrum of $(A_1,\dots,A_{d-1})$ acting on
$\L2{\B}$ iff $\chi$ is unitary. This can be used to give a \rep{}
theoretic definition of being Ramanujan, as in \Pref{S=sA}.
\begin{proof}[Proof of \Pref{S=sA}]
Assume every irreducible spherical infinite-dimensional \subrep{}
of $\L2{\dom{\Gamma}{G(F)}}$ is tempered. As $\Gamma$ is
cocompact, $\L2{\dom{\Gamma}{G(F)}}$ is a direct sum of
irreducible representations.  Let $f \in
\L2{\dc{\Gamma}{G(F)}{K}}$ be a non-trivial simultaneous
eigenvector of the $A_k$, with $A_k f = \lam_k f$.  By
\Pref{trade}, the $\lam_k$ are determined by some unramified
character $\chi$ . Consider $f$ as a $K$-fixed vector in
$\L2{\dom{\Gamma}{G(F)}}$. Since the only finite-dimensional
\rep{s} of $G(F)$ are the trivial ones, the \rep{} $\sg{f G(F)}$
is infinite-dimensional. Let $V$ be an irreducible quotient of
$\sg{f G(F)}$ in which $f \neq 0$, then $V$ is an irreducible
infinite-dimensional spherical \subrep{} of
$\L2{\dom{\Gamma}{G(F)}}$ so by assumption $V$ is tempered. It
then follows from \Pref{t=S} that $\chi$ is unitary, and so
$(\lam_1,\dots,\lam_{d-1}) \in \Aspec_d$.

In the other direction, let $V$ be an irreducible spherical
infinite-dimensional \subrep{} of $\L2{\dom{\Gamma}{G(F)}}$; then its
unique $K$-fixed vector $f$ is a simultaneous eigenvector of the
$A_k$, where $A_k f = \lam_k f$.  By assumption
$(\lam_1,\dots,\lam_d) \in \Aspec_d$.  The eigenvalues induce a
homomorphism $\omega = \omega_\chi$ for some unitary character
$\chi$, and by \Pref{sphrep}, $V$ is isomorphic
to a subquotient of $I_\chi$.  Consequently, $V$ is tempered.
\end{proof}

\subsection{Bounds on the spectrum of $A_k$ on $\L2{\dom{\Gamma}{\B}}$}\label{ss:bounds}

In the previous subsection we computed the spectrum of
$(A_1,\dots,A_{d-1})$ in their action on $\L2{\B}$. For a discrete
subgroup $\Gamma \leq G$, let
$\spec_{\dom{\Gamma}{\B}}(A_1,\dots,A_{d-1})$ denote the spectrum
of these operators in their action on $\L2{\dom{\Gamma}{\B}}$
(which is a finite set).

In this subsection we apply the classification of unitary \rep{s}
of $\GL[d](F)$ to give an upper bound on
$\spec_{\dom{\Gamma}{\B}}(A_1,\dots,A_{d-1})$ (which is
independent of $\Gamma$). In addition we state an Alon-Boppana
type theorem, due to W.~Li, that for suitable families of
quotients $\set{\dom{\Gamma_i}{\B}}$ of $\B$, $\overline{\cup
\spec_{\dom{\Gamma_i}{\B}}(A_1,\dots,A_{d-1})} \supseteq
\Aspec_d$.

Let $f \in \L2{\dc{\Gamma}{G}{K}}$ be a simultaneous eigenvector
of the $A_k$. Lift $f$ to $\L{2}{\dom{\Gamma}{G}}$, and recall
that the representation $\sg{fG}$ is unitary (since
the action of $G$ on $\L2{\dom{\Gamma}{G}}$ is unitary) and
spherical (since
$f$ is $K$-fixed).

The unitary
spherical \rep{s} were described by Tadic \cite{Tadic}, as
part of the classification of all the unitary
\rep{s} of $\GL(F)$. Such a spherical \rep{} is induced by a
character $\chi = \chi_1 \oplus \dots \oplus \chi_d$, where the
$\chi_i$ are combined into blocks. For the Satake parameters
$(z_{i_1},\dots,z_{i_s})$ of each block
$\chi_{i_1},\dots,\chi_{i_s}$, one of the following three options
holds: either $s = 1$ and $z_{i_1} \in S^1 = \set{z \in \C
\suchthat \abs[]{z} = 1}$;
$(z_{i_1},\dots,z_{i_s})$ is of the form
$$(q^{(s-1)/2}z,\dots,q^{(1-s)/2}z)$$
for $z \in S^1$; or (if $s = 2s'$ is even) it is of the form
$$(q^{(s'-1)/2+\alpha}z,\dots,q^{(1-s')/2+\alpha}z,q^{(s'-1)/2-\alpha}z,\dots,q^{(1-s')/2-\alpha}z)$$
for $z \in S^1$ and $0 < \alpha < 1/2$.

This set of possible parameters $(z_1,\dots,z_d)$ determines the
eigenvalues $(\lam_1,\dots,\lam_{d-1})$ via \Pref{fchi}. In
particular if $d \geq 3$, we obtain for the non-trivial
eigenvalues
$$\abs[]{\lam_k} \leq q^{k(d-k)/2} \cdot \s_k(q^{(d-2)/2},\dots,q^{(2-d)/2},1) \approx q^{k(d-k-\frac{1}{2})}$$
for every $k \leq d/2$ (and $\lam_{d-k} = \bar{\lam_k}$).

One can see that if $d \geq 3$ then $\abs[]{\lam_k} <
\binomq{d}{k}{q}$ for every non-trivial unitary representation
(where $\binomq{d}{k}{q}$ denotes the number of subspaces of
dimension $k$ in $\F_q^d$, which is the number of neighbors
of color $k$ of each vertex). In
particular the non-trivial eigenvalues of $A = A_1+\dots+A_{d-1}$
are bounded away from the trivial one. This demonstrates the fact
that $\PGL[d](F)$ has Kazhdan property $(T)$ and the quotient
graphs $\dom{\Gamma}{\B^1}$ are expanders for every $\Gamma$
\cite{Lu1}.

On the other hand, for $d = 2$ the eigenvalues
$q^{1/2}\s_1(q^{\alpha},q^{-\alpha})$ approach the degree $q+1$
when $\alpha \ra 1/2$, in accordance with the fact that
$\PGL[2](F)$ does not have property $(T)$.

For the lower bound, we quote
\begin{thm}[{\cite[Thm.~H]{Li2}}]
Let $X_i$ be a family of finite quotients of $\B$ with unbounded
injective radius. Then $\overline{\cup
\spec_{X_i}(A_1,\dots,A_{d-1})} \supseteq \Aspec_d$.
\end{thm}
This also follows from a multi-dimensional version of \cite{GZ}.

\subsection{Super-cuspidal and square-integrable \rep{s}}\label{ss:scsi}

Let $G$ denote the group $\GL(F)$ or $\PGL(F)$, and $Z = \Cent(G)$
its center. Let $\rho \co G \ra \End(V)$ be a unitary
representation. Recall that the
\defin{matrix coefficients} of $\rho$ are the functions
$\rho_{v,w} \co g \mapsto \innprod{\rho(g)v}{w}$ where
$v,w \in V$.
A unitary representation of $G$ is called
\defin{super-cuspidal}, if
its matrix coefficients are compactly supported modulo the center.
Notice
that the irreducible \rep{s} of $\GL[1](F)$ are all
super-cuspidal (as the group equals its center).

We say that a unitary \rep{} $\rho$ is
\defin{square-integrable}, if $\rho_{v,w} \in \L2{G/Z}$ for every
$v,w \in V$.
A \rep{} is square-integrable iff it is isomorphic to a \subrep{}
of $\L2{G}$ \cite[Prop.~9.6]{Knapp}.  Note that super-cuspidal
\rep{s} are square-integrable, and square-integrable \rep{s} are
tempered.

Let $s \divides d$ be any divisor, and let $P_s(F)$ denote the
parabolic subgroup corresponding to the partition of $d$ into $s$
equal parts. For a \rep{} $\psi$ of $\GL[d/s](F)$, we denote
\begin{equation}\label{Msdef}
M_s(\psi) = \Ind_{P_s(F)}^{\GL[d](F)} (\dd{F}{\frac{1-s}{2}}{\psi}
\,\HTplus\, \dd{F}{\frac{3-s}{2}}{\psi} \,\HTplus\, \cdots
\,\HTplus\, \dd{F}{\frac{s-1}{2}}{\psi}).
\end{equation}

The unique irreducible \subrep{} of $M_s(\dd{F}{(s-1)/2}{\psi})$
will be denoted by $C_s(\psi)$. It is known that if $\psi$ is
irreducible and super-cuspidal, then the induced \rep{}
$M_s(\dd{F}{(s-1)/2}{\psi})$ has precisely $2^{s-1}$ irreducible
subquotients, two of which (if
$s>1$) are unitary
\cite[p.~32]{HT} (notice that $M_1(\psi) =
\psi$). These subquotients are $\Cp_s(\psi)$, and a certain
irreducible quotient,
called the \defin{generalized Steinberg \rep{}} (or sometimes
"special \rep{}") and denoted by $\Sp_s(\psi)$.

\begin{prop}[{\cite[p.~32]{HT}, \cite{Zelevinsky}}]\label{sirep}
For $s > 1$ and $\psi$ an irreducible super-cuspidal \rep{} of
$\GL[d/s](F)$, $\Sp_s(\psi)$ is square-integrable, and
$\Cp_s(\psi)$ is not tempered.

Every square-integrable \rep{} of $\GL[d](F)$ is either
super-cuspidal, or of the form $\Sp_s(\psi)$ for a unique divisor
$s$ of $d$ and a unique super-cuspidal \rep{} $\psi$ of $\GL[d/s](F)$.
\end{prop}

\begin{rem}\label{moresirep}
If $s > 1$, $\Cp_s(\psi)$ is not tempered for any unitary
representation $\psi$.
\end{rem}

\begin{exmpl}\label{Cis1}
Let $\phi \co \mul{F} \ra \mul{\C}$ be a character, and $\psi =
\absdot[F]^{(1-d)/2} \phi$. Then $\Cp_d(\psi) = \phi\circ\det$,
which is one-dimensional.
\end{exmpl}
\begin{proof}
Let $B(F)$ denote the standard Borel subgroup of $\GL[d](F)$. By
definition, $\Cp_d(\psi)$ is the unique irreducible \subrep{} of
$$M_d(\phi) = \Ind_{B(F)}^{\GL[d](F)}( {{\absdot[F]^{(1-d)/2}\phi \,\HTplus\,
\cdots \,\HTplus\, \absdot[F]^{({d}-1)/2}}\phi}),$$
which is the unitary induction of $\Delta^{-1/2}\cdot
(\phi\circ\det)$ to $\GL(F)$. In particular, this \rep{}, when
restricted to $B(F)$, contains the \rep{} $\phi\circ\det$, which
is thus a \subrep{} of $M_d(\phi)$, so by
definition $\Cp_d(\psi)
= \phi\circ\det$.
\end{proof}

\Section{From local to global}\label{Prelim2}

\subsection{The global field}\label{ss:R0RT}

Let $k$ be a global field, $\Vals = \set{\nu}$ its nonarchimedean
discrete valuations, and $\Vals_{\infty}$ the Archimedean
valuations. For $\nu \in \Vals$, $k_{\nu}$ is the completion,
$\O_\nu = \set{x \suchthat \nu(x)\geq 0}$ the valuation ring of
$k_{\nu}$ (which is the closed unit ball of $k_{\nu}$ and thus
compact), and $P_{\nu} = \set{x \suchthat \nu(x)>0}$ the valuation
ideal.

Note that the ring of $\nu$-adic integers $k \cap \O_{\nu}$ of $k$
is a local ring, with maximal ideal $k \cap P_{\nu}$. Fix a
valuation $\nu_0 \in \Vals$, and set $F = k_{\nu_0}$. Consider the
intersection
$$
R_0 = \set{x\in k \suchthat \forall (\nu \in \Vals - \set{\valF})\ \nu(x)\geq 0} = \bigcap_{\nu \in \Vals \minusset
\set{\valF}} (k \cap \O_\nu).
$$

Recall that the valuations of $k = \F_q(y)$ are all
nonarchimedean. They are indexed by the prime polynomials of
$\F_q[y]$ and $1/y$. For a prime $p$ the valuation is $\nu_p(p^i
f/g) = i$ when $f$ and $g$ are prime to $p$, and the valuation
corresponding to $1/y$ is the minus degree valuation, defined by
$\yy{\nu}(f/g) = \deg(g)-\deg(f)$.
If $\nu_0 = \yy{\nu}$ then $R_0 = \F_q[y]$.

For every $x \in \mul{k}$ we have that
\begin{equation}\label{sumval}
\yy{\nu}(x) + \sum_{p}{\deg(p) \nu_p(x)} = 0,
\end{equation}
so $\valF(x) \leq 0$ for every $x \in R_0$. As a result $R_0$ is
discrete in $F$. It also follows that $R_0 \cap \O_{\nu_0} =
\F_q$.
For every $\nu$, choose a uniformizer $\uni_{\nu} \in R_0$, so
that $\nu(\uni_{\nu}) = 1$. Then, the completion of $k$ at $\nu$
is $\F_q\db{\uni_\nu}$, and the local ring of integers is
$\F_q[[\uni_\nu]]$. Note that if $\valF = \yy{\nu}$, we can choose
the uniformizers, $\uni_{\nu}$, to be the prime polynomials of
$\F_q[y]$, and $\uni_{\nu_0} = 1/y$.

Let $\I$ denote the set of functions $\ii \co \Vals \ra
\N\cup\set{0}$, such that $i_\nu = 0$ for almost all $\nu$ and
$i_{\nu_0} = 0$. The ideals of $R_0$ are indexed by functions $\ii
\in \I$, in the following way:
For $\ii \in \I$, we define
\begin{equation}\label{idealI}
I_{\ii} = \set{x \in R_0 \suchthat \nu(x) \geq i_{\nu}} =
\bigcap_{\nu \in \Vals \minusset \set{\nu_0}}(k \cap P_{\nu}^{i_{\nu}}),
\end{equation}
where we make the notational convention that ${P_\nu}^{0} =
\O_\nu$. {}From our choice of the uniformizers, it follows that if
$\valF = \yy{\nu}$, then $I_{\ii}$ is the (principal) ideal
generated by $\prod_{\nu \neq \nu_0}{\uni_{\nu}^{i_{\nu}}}$.
Notice that for the zero vector $\ii = 0$ ($i_{\nu} = 0$ for all
$\nu$), we obtain the trivial ideal $I_0 = R_0$.

Let $\times{k_\nu}$ be the direct product of the fields $k_\nu$
over all the valuations $\nu \in \Vals \cup \Vals_{\infty}$ of
$k$, and recall that the ring $\A$ of ad\`{e}les over $k$ is
defined
to be the restricted product
\begin{equation}\label{Adef}
\A = \set{x = (x_\nu) \in \times{k_{\nu}} \suchthat \nu(x_\nu)\geq 0 \quad
\mbox{for almost all $\nu$}}.
\end{equation}

The field $k$ embeds in $\A$ diagonally.
In a similar manner to the construction of $R_0$, we define
\begin{eqnarray}
\tilde{R}_0 & = & \set{(x_\nu) \in \A \suchthat \forall (\nu \in \Vals - \set{\valF})\ \nu(x_\nu)\geq 0} \label{R0t} \\
    & = & F \times \prod_{\nu \in \Vals -
\set{\valF}} \O_\nu \times \prod_{\nu \in \Vals_{\infty}} k_{\nu}.\nonumber
\end{eqnarray}
The ideals of finite index of $\tilde{R}_0$ are again indexed by
$\I$, and are of the form
\begin{equation}\label{idealIt}
\tilde{I}_{\ii} = \set{(x_{\nu}) \in \A \suchthat \forall(\nu \in \Vals - \set{\valF})\ \nu(x_\nu) \geq i_{\nu}},
\end{equation}
and with respect to the diagonal embedding, we have $R_0 = k \cap
\tilde{R}_0$ and $I_{\ii} = k \cap \tilde{I}_{\ii}$ for every $\ii
\in \I$. In fact, $\tilde{R}_0$ and $\tilde{I}_{\ii}$ are the topological
closures of $R_0$ and $I_{\ii}$, respectively.

\subsection{Strong Approximation}\label{ss:SA}

Let $\GG$
be a connected, simply connected, almost simple
linear
algebraic group, defined over $k$ (\eg\ $\SL$),
with a fixed
embedding into $\GL[r]$ for some $r$. For a subring $R$ of a
$k$-algebra $A$, we denote $\GG(R) = \GG(A) \cap \GL[r](R)$.
For
simplicity of notation (and as our applications are mainly for
positive characteristic), we assume $\GG(k_\nu)$ is compact for
all Archimedean places $\nu$.

The diagonal embedding $k \hra \A$, which is obviously discrete,
induces a discrete embedding $\GG(k) \hra \GG(\A)$.
Let $T$ be the set of valuations $\theta$ such that
$\GG(k_{\theta})$ is compact; this is a finite set \cite{PR}. Fix
a valuation $\valF \in \Vals - T$, and let $F = k_{\nu_0}$ denote
the completion with respect to this special valuation. $\GG(k)$ is
a lattice of finite co-volume in $\GG(\A)$, and moreover if $T
\neq \emptyset$, $\GG(k)$ is a cocompact lattice
\cite[Thm.~5.5]{PR}.

\begin{thm}[{Strong Approximation \cite{Prasad}, \cite{PR}}]
The product $\GG(k) \GG(F)$ is dense in $\GG(\A)$.
\end{thm}

So for every open subgroup $U$ of $\GG(\A)$,
\begin{equation}\label{dense}
\GG(k)\GG(F) U = \GG(\A).
\end{equation}

\begin{cor}\label{SAPGL}
Let $U \sub \GG(\A)$ be a compact subgroup such that $\GG(F) U$ is
open,
and $\GG(F) \cap U = 1$. Set $\GammaU = \GG(k) \cap \GG(F)U.$ Then
its projection to $\GG(F)$ (which we will also denote by
$\GammaU$) is discrete, and
\begin{equation}\label{denseA}
\dc{\GG(k)}{\GG(\A)}{U} \isom \dom{\GammaU}{\GG(F)}.
\end{equation}
\end{cor}

For example, if $U = \prod_{\nu \in \Vals-\set{\valF}} \GG(\O_{\nu}) \times
\prod_{\nu \in \Vals_{\infty}} \GG(k_{\nu})$,
then $\Gamma = \GammaU$ is the arithmetic subgroup $\GG(R_0)$. More
generally,
let $\ii = (i_\nu) \in \I$ be a function
corresponding to an ideal $\tilde{I}_\ii$, and let
$$U_{\ii} = \prod_{\nu \in \Vals - \set{\valF}} \con{\GG}{\O_{\nu}}{P_{\nu}^{i_{\nu}}} \times \prod_{\nu \in \Vals_{\infty}} \GG(k_{\nu})$$
where $\con{\GG}{\O_{\nu}}{P_{\nu}^{i_{\nu}}} =
\Ker(\GG(\O_{\nu}) \ra \GG(\O_{\nu}/P_{\nu}^{i_{\nu}}))$ is a congruence
subgroup.
Then $\GG(F) U_{\ii} = \con{\GG}{\tilde{R_0}}{\tilde{I}_{\ii}} =
\Ker(\GG(\tilde{R_0})) \ra \GG(\tilde{R_0}/\tilde{I}_{\ii})$
is an open subgroup of $\GG(\A)$, and we set
\begin{equation}\label{Gammadef}
\Gamma_{\ii} = \con{\GG}{R_0}{I_{\ii}} = \GG(k) \cap \GG(F)U_\ii,
\end{equation}
called the \defin{principal congruence subgroup mod $I_\ii$} of $\GG(R_0)$.
Again, when $T \neq \emptyset$, this is a cocompact
lattice
in $\GG(F)$.

\subsection{Automorphic representations}\label{sss:localglobal}

The group $\GG(\A)$ acts on the space $\automrep[\GG]$ by
multiplication from the right. The sub-modules are called
automorphic representations of $\GG(\A)$. The closed irreducible
sub-modules are said to be discrete, or to belong to the discrete
spectrum
. Its complement is called the continuous spectrum.  If $T \not =
\emptyset$ then there is no continuous spectrum.

Let $K_\nu = \GG(\O_{\nu})$
and recall that for every $(g_{\nu}) \in \GG(\A)$, $g_{\nu} \in
K_{\nu}$ for almost all $\nu$.
Given irreducible representations $\pi_{\nu} \co \GG(k_{\nu}) \ra
\End(V_{\nu})$, with all but finitely many being
$K_{\nu}$-spherical, one defines the restricted tensor product
$\pi = \tensor[] \pi_{\nu} \co \GG(\A) \ra \End(\tensor[]'
V_{\nu})$ \cite{Bump}.

A fundamental theorem \cite[Thm.~3.3.3]{Bump}
states that any irreducible
automorphic \rep{} of $\GG(\A)$ is isomorphic to such a restricted
tensor product. The representations $\pi_\nu$ in $\pi =
\tensor\pi_{\nu}$ are called the (local) components of $\pi$, and
since $\pi$ is irreducible, they are also irreducible. Moreover,
$\pi$ is admissible iff all its components are.

\begin{prop}\label{trivtriv}
Assume that $\GG(k_{\nu})$ is non-compact. If the
component $\pi_{\nu}$ of an irreducible automorphic
\rep{} $\pi$ of $\GG(\A)$ is trivial, then $\pi$ is trivial.
\end{prop}
\begin{proof}
Let $\pi = \tensor\pi_{\nu}$ be an automorphic representation
acting on $V \leq \automrep[\GG]$, where $\pi_{\nu}$ is trivial.
For $f \in V$, (assumed to be $K$-finite, see
\cite[Thm.~3.3.4]{Bump}),
$f$ is $\GG(k_{\nu})$-invariant from the right and
$\GG(k)$-invariant from the left.  However by Strong
Approximation, $\GG(k_{\nu})\GG(k)$ is dense, so $f$ must be
constant everywhere, making $\pi$ trivial.
\end{proof}

Recall that an irreducible \rep{} is $U_{\ii}$-spherical if it has a
$U_{\ii}$-fixed vector.
We assume $\GG(F)$ is non-compact where $F = k_{\valF}$.

\begin{prop}\label{locallift}
Let $\pi$ be an irreducible, $U_{\ii}$--spherical automorphic
\rep{} of $\GG(\A)$. Then $\pi_{\nu_0}$ is a \subrep{} of
$\L2{\dom{\Gamma_{\ii}}{\GG(F)}}$.

Conversely, if $\rho \leq \L2{\dom{\Gamma_{\ii}}{\GG(F)}}$ is
irreducible, then there exists an irreducible $U_{\ii}$--spherical
automorphic \rep{} $\pi$ of $\GG(\A)$ such that $\pi_{\nu_0}$ is
isomorphic to $\rho$.
\end{prop}

The second assertion is seen by lifting a function $f \in V_\rho$
(where $V_{\rho}$ is the \rep{} space) from
$\dom{\Gamma_{\ii}}{\GG(F)}$ to $\dom{\GG(k)}{\GG(\A)}$ using
\Cref{SAPGL}, and taking $\pi$ to be an irreducible quotient of
the (right) $\GG(\A)$-module generated by $f$.

\medskip

\subsection{The conductor}\label{sss:conductor}

For a \rep{} $\rho$ of $\GG(k_{\nu})$, the \defin{conductor} of
$\rho$, $\cond{\rho} = i$, is defined to be
the minimal $i \geq 0$, for which there is a
$\con{\GG}{\O_{\nu}}{P_{\nu}^i}$-fixed vector in $V$
(such an $i$ exists since the \rep{} is admissible).
In
particular, $\cond{\rho} = 0$ iff $\rho$ is spherical.

Now let $\pi$ be an irreducible automorphic representation of
$\GG(\A)$.
Since almost all the local components are spherical,
$\cond{\pi_{\nu}} = 0$ for almost every $\nu$. We thus let
$\cond{\pi}$ be the function $\ii \co \Vals \ra \N\cup \set{0}$
defined by $i_{\nu} =\cond{\pi_{\nu}}$ (note that
$\ii$ is not in $\I$ in general, as we do not assume $i_{\nu_0} = 0$).
\begin{rem}\label{cond=U}
Let $\ii = \cond{\pi}$. Then $H =
\con{\GG}{\O_{\valF}}{P_{\valF}^{i_{\valF}}}U_{\ii}$ is the
maximal principal congruence subgroup for which $\pi$ has an
$H$-fixed vector.
\end{rem}

\medskip

The results of this section will be used later for non-simply
connected cases, which requires some minor modifications.  Let $G$
be a connected, almost simple algebraic group over $k$.  Let $\Pi
\co \GG \ra G$ be its simply connected cover.  Then $\Pi(\GG(\A))
\normal G(\A)$ and the quotient is abelian (of finite exponent).
In this situation, \Pref{trivtriv} becomes

\begin{prop}\label{trivtriv2}
Assume that $G(k_{\nu})$ is non-compact. If the component
$\pi_{\nu}$ of an irreducible automorphic \rep{} $\pi$ of $G(\A)$
is one dimensional, then $\pi$ is one dimensional.
\end{prop}

\Section{Global automorphic \rep{s}}\label{Prelim3}

Let $G$ be an almost simple, connected algebraic group defined
over $k$, where $k$ is a global field of arbitrary characteristic.
The discrete spectrum of automorphic \rep{s} is
composed of cuspidal and residual \rep{s}. The cuspidal
representation space is comprised of functions $f \in \L2{\dom{G(k)}
{G(\A)}}$ which satisfy
$\int_{\dom{N(k)}{N(\A)}} f(ng)dn = 0$ for every $g \in G(\A)$ and
for every $N$, where $N$ is a unipotent radical of a parabolic
subgroup of $G$. Since the cuspidal condition involves integration
from the left, and the action is by right translation, this is a
sub-representation space. The other discrete irreducible \rep{s}
are called residual.

Recently, L.~Lafforgue has proved the following version of the
Ramanujan conjecture:
\begin{thm}[\cite{Lafforgue}, \cite{Rev}] \label{Lafforgue}
Assume $k$ is of positive characteristic and $G=\GL[d]$.  Let $\pi
= \otimes \pi_\nu$ be an irreducible, cuspidal representation with
finite central character.  For all $\nu$, if $\pi_\nu$ is
spherical then $\pi_\nu$
is tempered.
\end{thm}

\subsection{The residual spectrum}\label{residual}

All the one-dimensional \rep{s} are residual, and when $G=\GL$
and $d$ is prime these are the only residual \rep{s}. If $d$ is
not a prime, the other residual
\rep{s} can be described in terms of the cuspidal representations
of smaller rank, as follows:

An element $(a_\nu) \in \A$ is invertible only if for almost all
$\nu$, $a_\nu \in \mul{\O_{\nu}}$. We can thus define an absolute
value on $\mul{\A}$ by $\abs[\A]{(a_{\nu})} = \prod \abs[k_{\nu}]{a_\nu}$,
which is a finite product. The modular function for parabolic
subgroups of $\GL(\A)$ is defined as in the local case (see
\Eq{modularP}, with $\dd[a]{\A}{}$ for each block),
and likewise we have a unitary induction from
parabolic subgroups, with similar properties to the local case.

Let $s > 1$ be a divisor of $d$, and let $\pi$ be any cuspidal
automorphic \rep{} of $\GL[d/s](\A)$. The representation
\begin{equation}\label{TsdA}
T_{s}(\pi) =
\Ind_{P_s(\A)}^{\GL(\A)}(
\dd{\A}{\frac{1-s}{2}}{\pi} \oplus \dd{\A}{\frac{3-s}{2}}{\pi} \oplus \cdots \oplus \dd{\A}{\frac{s-1}{2}}{\pi}
)
\end{equation}
has a unique irreducible \subrep{} $J(T_{s}(\pi))$
(here $P_s(\A)$ is the parabolic subgroup of $\GL(\A)$
associated to the decomposition into $s$ blocks of size $d/s$).

\begin{thm}[\cite{MW}]\label{MW}
The residual spectrum of $\L2{\dom{\GL(k)}{\GL(\A)}}$ consists of
the representations $J(T_{s}(\pi))$ for proper divisors $s
\divides d$ and $\pi$ a cuspidal \rep{} of $\GL[d/s](\A)$.
\end{thm}

Comparing Equations \eq{Msdef} and \eq{TsdA}, the local
$\nu$-component of  $J(T_s(\pi))$
is seen to be the (unique) irreducible \subrep{} of
$M_s(\pi_{\nu})$,
namely $\Cp(\pi_{\nu})$ which was defined in Subsection \ref{ss:scsi}.
{}From \Rref{moresirep} we then obtain
\begin{cor}\label{Laff}
(a) Every local component of a residual representation is
non-tempered.

(b) If $\pi$ is an irreducible automorphic representation of $GL_d$
where one of its local components is tempered, then $\pi$ is
cuspidal (and in positive characteristic, all of its spherical components are tempered by \Tref{Lafforgue}).
\end{cor}

\subsection{The \JLc}\label{ss:JL}

Let $D$ be a division algebra of degree $d$ over $k$, and let
$D_{\nu} = D \tensor[k] k_{\nu}$. Then by the
Albert-Brauer-Hasse-Noether theorem, $D_{\nu} \isom
\M[d](k_{\nu})$ for almost every completion $k_{\nu}$. Let $G' =
\mul{D}$, which is a form of inner type of $G = \GL[d]$.
Let $T$ denote the (finite) set of
valuations $\theta$ such that $D \tensor k_{\theta}$ is not split.
We assume that for every $\theta \in T$, $D \tensor k_{\theta}$ is
a division algebra.
There is an injective correspondence, called the local \JLc, which
maps every irreducible, unitary representation $\rho'$ of
$G'(k_{\theta})$ ($\theta \in T$) to an irreducible, unitary
square-integrable (modulo the center) representation $\rho =
\JL_\theta(\rho')$ of $G(k_{\theta})$ (see \cite{Rogawski} or
{\cite[p.~29]{HT}} for details).

If $\phi$ is a character of $\mul{k_{\theta}}$, then \cite[p.~32]{HT}
\begin{equation}\label{JLchar}
\JL_{\theta}(\phi \circ \det) = \Sp_d(\absdot[k_{\theta}]^{(1-d)/2}\phi)
\end{equation}
where $\Sp_d$ is defined in Subsection \ref{ss:scsi}.
Recall by \Eref{Cis1}, that
$\Cp_d(\absdot[k_{\theta}]^{(1-d)/2} \phi)$ is a one-dimensional \rep{}.

The global \JLc\ maps an irreducible automorphic representation
$\pi'$ of $G'(\A)$ to an irreducible automorphic representation
$\pi = \JL(\pi')$ of $G(\A)$ which occurs in the discrete spectrum
(see \cite[p.~195]{HT}).
If $\nu \not \in T$, then
$$\JL(\pi')_\nu \isom \pi'_\nu.$$
Note that  the restrictions of $\cond{\pi}$ and $\cond{\pi'}$ to
$\Vals-T$ are equal.

The situation in the other local components is as follows: let
$\theta \in T$, and consider the component $\pi'_{\theta}$ of
$\pi'$. The local \JLc\ maps $\pi'_{\theta}$ to an irreducible
square-integrable \rep{} $\JL_{\theta}(\pi'_{\theta})$ of
$G(k_{\theta})$, which is by \Pref{sirep} a generalized Steinberg
\rep{,} of the form $\Sp_s(\psi)$ for some divisor $s \divides d$
and super-cuspidal \rep{} $\psi$ of $\GL[d/s](k_{\theta})$. Then
$\JL(\pi')_{\theta}$ is isomorphic to either $\Sp_s(\psi)$ or
$\Cp_s(\psi)$.

\begin{thm}[{\cite[p.~196]{HT}}]\label{HTJL}
The image of $\JL$ (for a fixed $D$) is the set of irreducible
automorphic representations $\pi$ of $\GL(\A)$ such that $\pi$
occurs in the discrete spectrum and for every $\theta \in T$ there
is a positive integer $s_\theta | d$ and an irreducible
super-cuspidal representation $\psi_\theta$ of
$\GL[d/s_\theta](k_\theta)$ such that $\pi_\theta$ is isomorphic
to either $\Sp_{s_\theta}(\psi_\theta)$ or
$\Cp_{s_\theta}(\psi_\theta)$.
\end{thm}

Throughout the  book \cite{HT}, the authors assume characteristic
zero. However, see \Rref{valid}.

\Section{Proofs of Theorems \ref{IntroT1} and
\ref{IntroT2}}\label{inner}

Let $k$ be a global field of prime characteristic, $D$ a division
algebra of degree $d$ over $k$, $G'$ the algebraic group
$\mul{D}/\mul{Z}$ where $Z$ is the center, and $G = \PGL[d]$.

Let $T$ denote the set of ramified primes, namely valuations
$\theta$ for which $D_{\theta}  = D \tensor k_\theta$ is
non-split. We again
assume that for such primes, $D_{\theta}$ is a division algebra.
It follows that $G'(k_\theta)$ is compact for $\theta \in T$. The
valuation $\theta$ extends uniquely to a valuation of
$D_{\theta}$, and we let $\O_{D_{\theta}}$ denote the ring of
integers there.

The group $G'(\O_{\theta})$ depends on the specific embedding
$G'(k) \hra \GL[r](k)$, namely, $G'(k_{\theta})$ is the subgroup
of $\GL[r](k_{\theta})$ defined by the equations defining $G'(k)$,
and $G'(\O_{\theta}) = G'(k_{\theta}) \cap \GL[r](\O_{\theta})$.
For most of our applications the precise embedding is irrelevant
($G'(\O_{\theta})$ is well defined up to commensurability anyway).
However, for \Tref{IntroT2}(b), we need the embedding to satisfy
\begin{equation}\label{GOt}
G'(\O_{\theta}) \sup
\mul{k_{\theta}}\mul{\O_{{D_{\theta}}}}/\mul{k_{\theta}},
\end{equation}
where both groups are viewed as subgroups of $G'(k_{\theta}) =
\mul{(D\tensor[k] k_{\theta})}/\mul{k_{\theta}}$, which is
embedded in $\GL[r](k_{\theta})$ for some $r$.
This condition is in fact satisfied by a natural embedding. Let
$E$ be a cyclic extension of dimension $d$ over $k$, which is
unramified at every $\theta \in T$ (the existence of $E$ is
guaranteed by Grunwald's theorem for function fields
\cite[Chap.~10]{AT}).
{}From Albert-Brauer-Hasse-Noether theorem it follows that $E$ is
a splitting field of $D$, making $D$ a cyclic division algebra.
Moreover there is an element $z \in D$ such that $D = E[z]$ and
conjugation by $z$ is an automorphism of $E$, generating
$\Gal(E/k)$.

Let $e_1,\dots,e_d$ be an integral basis of $E/k$ (with respect to
every $\theta \in T$). Then, for every valuation $\theta \in T$,
$\O_{E_{\theta}} = \sum \O_{\theta} e_i$, where $E_{\theta} = E
\tensor[k] k_{\theta}$ and $\O_{E_{\theta}}$ is its ring of
integers. Now, $z$ can be chosen so that $\O_{D_\theta} =
\O_{E_{\theta}}[z] = \sum_{i,j} \O_{\theta}e_iz^j$. The left
regular representation of $D$ via the basis $\set{e_i z^j}$
defines an embedding $\mul{D} \ra \GL[d^2](k)$ which sends
$\mul{\O_{D_{\theta}}}$ into $\GL[d^2](\O_{\theta})$ (and central
elements to scalar matrices). Composing this with the adjoint
representation of $\PGL[d^2](k)$ (into $\GL[d^4](k)$), we obtain
an embedding of $G'(k) = \mul{D}/\mul{k}$ which satisfies
\Eq{GOt}.

\begin{lem}\label{GDGO}
For $\theta \in T$,
$\mul{\O_{D_{\theta}}}$ is normal in $\mul{D_{{\theta}}}$, and
$\mul{D_{{\theta}}}/\mul{k_{\theta}}\mul{\O_{D_{\theta}}} \isom \Z/d$.
\end{lem}
\begin{proof}
The uniformizer $\uni$
of $k_{\theta}$ is a uniformizer for $E_{\theta}$ as well,
and (by choosing the generator
$\s \in \Gal(E_{\theta}/k_{\theta})$ appropriately)
we may assume $D_{\theta} = E_{\theta}[z]$ where $z^d = \uni$ and conjugation
by $z$ induces $\s$. Since $z$
normalizes $\mul{\O_{{D_{\theta}}}}$,
this is a normal subgroup of $\mul{D_{{\theta}}}$.

The elements of value zero in $\O_{D_{\theta}}$ are invertible there, so
every element of $\mul{D_{{\theta}}}$ is
of the form $c z^i$ for some $c \in \mul{\O_{D_{\theta}}}$ and an integer $i$.
Such an element is equivalent to $z^i$ in
$\mul{D_{{\theta}}}/\mul{k_{\theta}}\mul{\O_{D_{\theta}}}$, and
$z^d = \uni \in \mul{k_{\theta}}$. Finally, $z^i$ induces a
non-trivial automorphism on $E$ for for every $0<i<d$, so the
order of $z$ modulo the center is equal to $d$.
\end{proof}

By our assumption \eq{GOt}, the lemma implies that
$G'(k_{\theta})/G'(\O_{\theta})$ is a quotient of $\Z/d$.

\medskip

For $\ii \in \I$ set $\Gamma_{\ii} = \con{G'}{R_0}{I_{\ii}}$,
as in \Eq{Gammadef}.
For $\ii,\jj \in \I$, we say that $\ii \leq_T \jj$ if $i_{\nu}
\leq \jj_{\nu}$ for every $\nu \in \Vals \minusset T$.

\begin{prop}\label{3.1}
Let $\ii \in \I$.
The complex $\dom{\Gamma_{\ii}}{\B}$
is Ramanujan
iff every spherical
infinite-dimensional $\nu_0$-component of an
irreducible automorphic discrete \rep{} $\pi'$ of $G'(\A)$ with
$\cond{\pi'} \leq_{T} \ii$, is tempered.
\end{prop}
\begin{proof}
This follows immediately from Propositions \ref{S=sA} and \ref{locallift}
(and Remark \ref{cond=U}).
\end{proof}

We can now prove the theorems stated in the Introduction.

\begin{proof}[Proof of Theorems \ref{IntroT1} and \ref{IntroT2}(a)]
Write the given ideal of $R_0$ as $I = I_{\ii}$ for $\ii \in \I$
(see \Eq{idealI}). Let $\pi' $ be an irreducible discrete
automorphic \rep{} of $G'(\A)$ with $\cond{\pi'} \leq_T \ii$, and
assume $\rho = \pi'_{\valF}$ is spherical and
infinite-dimensional. By \Pref{3.1}, $\dom{\Gamma(I)}{\B}$ is
Ramanujan iff in all such cases $\pi'_{\valF}$ is tempered.

By the \JLc, there is an irreducible automorphic \subrep{} $\pi$ of
$\L2{\dom{G(k)}{G(\A)}}$ such that $\pi_{\nu} = \pi'_{\nu}$ for
every $\nu \not \in T$. In particular, $\pi_{\valF} =
\pi'_{\valF}$.

Assume $d$ is prime, then all the infinite dimensional automorphic
\rep{s} of $G(\A)$ are cuspidal, so $\pi$ is cuspidal. By
Lafforgue's \Tref{Lafforgue}, the components of a cuspidal
representation are tempered. Therefore, $\rho = \pi_{v_0}$ is
tempered, and \Tref{IntroT1} is proved.

Now let $d$ be arbitrary, and assume $i_{\theta} = 0$ for some
$\theta \in T$ (namely $I$ is prime to $\theta$). Thus,
$\pi'_{\theta}$ has a $G'(\O_{\theta})$ - fixed vector. By Lemma
\ref{GDGO}, $G'(\O_{\theta})$ is normal in $G'(k_{\theta})$, and
$\pi'_{\theta}$ is an irreducible representation of the cyclic
quotient, so it is one-dimensional. Write $\pi'_{\theta} = \phi
\circ \det$ for a suitable character $\phi \co \mul{k_{\theta}}
\ra \C$ (of order $d$), where here $\det$ stands for the reduced
norm of $G'(k_{\theta})$.

By \Eq{JLchar}
we have that $\JL_{\theta}(\pi'_{\theta}) =
\JL_{\theta}(\phi \circ \det) = \Sp_d(\psi)$ for the character $\psi
= \absdot^{(1-d)/2} \phi$ of $\mul{k_{\theta}}$. By \Eref{Cis1},
$\Cp_d(\psi)$ is one-dimensional.

As mentioned in Subsection \ref{ss:JL}, $\pi_{\theta}$ is
isomorphic to either $\Sp_s(\psi)$ or $\Cp_s(\psi)$, but
$\pi_{\theta}$ cannot be one-dimensional (by Remark
\ref{trivtriv2}).  Therefore, $\pi_{\theta} = \Sp_s(\psi)$, which
is square-integrable (\Pref{sirep}) and, in particular, tempered.
Now, \Cref{Laff}(a) implies that $\pi$ is cuspidal, and by
\Tref{Lafforgue}, $\pi_{\valF} = \rho$ is tempered too.
\end{proof}

\begin{proof}[Proof of \Tref{IntroT2}(b)]
By \Pref{3.1}, we need to find an irreducible sub-representation $\pi'$ of
$\L2{\dom{G'(k)}{G'(\A)}}$ such that $\pi'_{\nu_0}$ is spherical
and non-tempered. Then $\pi'_{\nu_0}$ would be a \subrep{} of
$\L2{\dom{\con{G'}{R_0}{I}}{G(F)}}$ for some $I \normali R_0$, and
$\dom{\con{G'}{R_0}{I}}{\B}$ would not be Ramanujan.

We use the following result, which is a variant of a special case of
{\cite[Thm.~2.2]{Vig}}.
\begin{prop} \label{wVig}
Let $T = \set{\theta_1,\dots,\theta_t}$ and $\nu_1 \not \in T$ be
valuations of $k$.
For $i = 1,\dots,t$, let $\psi_i$ be a super-cuspidal \rep{} of
$\PGL[m](k_{\theta_i})$, where $m>1$ is fixed.

Then, there exists an automorphic cuspidal \rep{} $\;\pi\,$ of $\PGL[m](\A)$,
such that $\pi_{\theta_i} = \psi_i$ for
$i = 1,\dots,t$, and $\pi_{\nu'}$ is spherical for every valuation
$\nu' \not \in T \cup \set{\nu_1}$.
\end{prop}
\begin{proof}
Here we let $G$ denote the group $\PGL[m]$.
Let $f_{\theta_i}$ be matrix coefficients of $\psi_i$, and let
$U_{\theta_i}$ denote the (compact and open) support.
For $\nu \not \in T \cup \set{\nu_1}$
let $U_{\nu} = G(\O_{\nu})$, and choose an open compact
subgroup $U_{\nu_1}$ of $G(k_{\nu_1})$ such that
$U = \prod{U_{\nu}} \sub G(\A)$ intersects $G(k)$ only in the
identity.
For $\nu \neq \theta_1,\dots,\theta_t$, let $f_{\nu}$
be the characteristic function of $U_{\nu}$.
Let $f = \tensor f_{\nu} \in \L2{G(\A)}$.

Define an operator $R_f \co \L2{\dom{G(k)}{G(\A)}} \ra \L2{\dom{G(k)}{G(\A)}}$ by
$$R_f\varphi(g) = \int_{G(\A)}{f(g^{-1}x)\varphi(x)dx}.$$
The image of $R_f$ is in the discrete spectrum. Let $\pi$ be an
irreducible representation in the image, then $\pi_{\theta_i} = \psi_i$
and in particular $\pi$ is cuspidal. Moreover, $f_{\nu'}$ is a fixed vector
of $\pi_{\nu'}$ so these are spherical for
every $\nu' \not \in T \cup \set{\nu_1}$.
It remains to show that $R_f \neq 0$:
\begin{eqnarray*}
R_f\varphi(g)
    & = & \int_{\dom{G(k)}{G(\A)}} K_f(g,x) \varphi(x) dx
\end{eqnarray*}
where $K_f(g,x) = \sum_{\gamma \in G(k)} f(g^{-1}\gamma x)$, which is
a finite sum since $f$ is compactly supported. But
$K_f(1,1) = f(1) + \sum_{1\neq \gamma \in G(k)}f(\gamma) = 1$,
showing that $K_f \neq 0$ and $R_f \neq 0$.
\end{proof}

For $T$ we take the usual set of places in which $D$ remains a division
algebra, and we choose an arbitrary $\nu_1 \not \in T \cup \set{\valF}$.

Now pick any proper divisor $s$ of $d$.
For every $i = 1,\dots,t$ choose a super-cuspidal \rep{} $\psi_i$ of
$\PGL[d/s](k_{\theta_{i}})$,
and let
$\pi$ be the \rep{} of $\PGL[d/s](\A)$ given by \Pref{wVig};
in particular $\pi_{\valF}$ is
spherical. Then let $\tilde{\pi} = T_{s}(\pi)$, as in \Eq{TsdA},
and let $\opi = J(\tilde{\pi})$ be its
unique irreducible \subrep{}.
By \Pref{MW}, $\opi$ is in the residual spectrum, and
in particular $\opi_{\valF}$ is spherical (since $\valF \neq \nu_1$)
and non-tempered (\Cref{Laff}(a)).

Now, for every $i = 1,\dots, t$, $\opi_{\theta_i} =
\Cp_s(\dd{F}{(1-s)/2}\psi_i)$ (see the remark preceding
\Cref{Laff}),
so by \Tref{HTJL}, $\opi$ is in the image of the \JLc,
corresponding to a \rep{} $\opi'$ of $G'(\A)$ where $G' =
\mul{D}/\mul{Z}$. But $\opi_{\nu_0} = \opi'_{\nu_0}$, so this
component is spherical and not tempered.
\end{proof}

\Section{Outer forms}\label{outer}

Theorem \ref{IntroT1} (especially when compared to \Tref{IntroT2}(b))
may suggest that if $d$ is an odd
prime, then every finite quotient of
the Bruhat-Tits building $\B = \B_d(F)$ is Ramanujan,
where $F$ is a local field.

Indeed, if $Y$ is such a finite quotient of $\B$, then the
fundamental group $\Gamma_1 = \pi_1(Y)$ acts on $\B$, the
universal cover of $Y$, and $Y = \dom{\Gamma_1}{\B}$. By a well
known result of Tits, $\Aut(\B)$ is $G = \PGL[d](F)$, up to
compact extension.
It seems likely that $\Gamma_1$ has a subgroup of finite index
$\Gamma$ which is contained in $G$, and the corresponding finite
cover of $Y$ can be obtained as $\dc{\Gamma}{G}{K}$. Now, by
Margulis' arithmeticity theorem \cite{Mabook}, $\Gamma$ is an
arithmetic lattice of $G$.

A well known conjecture of Serre \cite{SerreCSP} asserts that
arithmetic lattices of $G$ (where $d \geq 3$) satisfy the
congruence subgroup property. This essentially means that every
finite index subgroup is a congruence subgroup. If $\Gamma$ is of
inner type, our \Tref{IntroT1} applies to it, and shows that the
quotients are really Ramanujan. However, there are other
arithmetic subgroups (see for example the classification of the
$k$-forms of $\GL[d]$ in \cite[III.1.4]{Serre}).

The outer forms of $\PGL[d]$
all come from the following general construction:
let $k$ be a global field, $k'/k$ a quadratic separable extension, and $A$
a $k'$-central simple algebra with an involution $u \mapsto u^*$
which induces the non-trivial automorphism of $k'/k$ on the center
of $A$.  Let $\Norm[k'/k]$ denote the norm map.  The algebraic
group $G' = \set{u \in A \suchthat uu^* = 1} / Z$ (where $Z =
\Ker(\Norm[k'/k])$ is the center) gives a form of $\PGL[d]$. Now,
if $d$ is a prime, $A$ may be either a division algebra, or the
matrix algebra $\M(k')$. The second case corresponds to Hermitian
forms \cite{PR}, \ie\ $G'$ is
$\PGU(q,k') = \set{a \in \M(k')\suchthat q(a(v)) = q(v)} / Z$ of
operators preserving the Hermitian form $q \co (k')^d \ra k'$. In
this situation, the involution on $A$ is $a \mapsto
b^{-1}\bar{a}^t b$, where $b$ is a skew-symmetric matrix
representing $q$.

But, if $\mychar{k} = p > 0$, every Hermitian form over $k$ represents $0$
if $d \geq 3$.
Indeed this is known to be true for local fields \cite[Sec.~4.2]{Schar}
and by Hasse Principal \cite[Sec.~4.5]{Schar}, this is also true for $k$.
Now in order to form a cocompact arithmetic lattice $\Gamma$ in
$\PGL(F)$, the form $G'$ should be anisotropic (\ie\ have $k$-rank
zero), but if $q$ represents $0$ over $k$, the $k$-rank is greater
than zero.
Thus, there are no arithmetic lattices of Hermitian form type if $d \geq 3$
and $F$ is of positive characteristic (the situation is different for
characteristic zero, see below).

On the other hand, the case when $A$ is a division algebra is possible
(\eg\ the cyclic algebra $A = \F_{q^d}(t)[z \subjectto z^d = t]$ where
$z$ induces the Frobenius automorphism on $\F_{q^d}$, is a division algebra
with center $k' = \F_q(t)$, and has an involution defined by $z^* = z^{-1}$
and $\alpha^* = \alpha$ for $\alpha \in \F_{q^d}$,
which is non-trivial on $k'$).
We do not know if \Tref{IntroT1} is valid in this case, but
if it is true then together with Serre's conjecture this would imply the
remarkable possibility that if $\mychar F > 0$ and $d \geq 3$ is a prime,
then all the finite quotients of $\B_d(F)$ are Ramanujan.
We leave it, however, as an open problem.

For $d = 2$, \ie\ $\PGL[2](F)$, all arithmetic lattices are of
inner type, as the Dynkin diagram of $A_1$ does not have graph
automorphisms, so Theorem \ref{IntroT1} applies for all lattices
(a result which has been proved before by Morgenstern
\cite{Morg}). Still we have

\begin{prop}\label{misref}
If $d = 2$, for every non-archimedean local field $F$, of any characteristic
$\PGL[2](F)$ has cocompact
(arithmetic) lattices,
such that the quotient  $\dom{\Gamma}{\B_2(F)}$ of the tree
$\B_2(F)$  is not Ramanujan.
\end{prop}
\begin{proof}
The group $\PGL[2](F)$ has cocompact (arithmetic) lattices, and
these are virtually free (\cf \cite{Serre:Trees}). Let $\Gamma$ be
a free cocompact lattice in $\PGL[2](F)$, so $\Gamma' =
[\Gamma,\Gamma]$ is of infinite index in $\Gamma$. Let $\Gamma_n$
be a sequence of finite index subgroups of $\Gamma$, such that
$\bigcap \Gamma_n = \Gamma'$. By \cite[Sec.~4.3]{alexbook}, the
graphs $\dom{\Gamma_n}{\B}$ are not expanders, let alone Ramanujan
graphs
. Of course, in light of \Tref{IntroT1} (or \cite{Morg}) for positive
characteristic, and \cite[Thm.~7.3.1]{alexbook} (see also \cite{JL})
for zero characteristic,
almost all the $\Gamma_n$ are non-congruence subgroups.
\end{proof}

\begin{lem}
Let $F$ be a local nonarchimedean field of characteristic zero.
For every $d \geq 2$, there exists a number field $k$ with a quadratic
extension $k'$ such that $k \sub k' \sub F$, and an anisotropic
Hermitian form $q$ of dimension $d$ over $k'/k$.
\end{lem}
\begin{proof}
Let $p$ be the prime such that $\Q_p \sub F$. Choose a natural number
$\delta > 0$ such that $-\delta$ is a quadratic residue modulo $p$ if $p$
is odd (\eg\ $\delta = p-1$), and take $\delta = 7$ if $p = 2$.
Let $k = \Q$
and $k' = \Q[\sqrt{-\delta}]$,
and let $u \mapsto \bar{u}$ denote the
non-trivial automorphism of $k'/k$.
Let $q(u_1,\dots,u_d) = u_1 \bar{u_1} + \dots + u_d \bar{u_d}$. Writing
$u_i = x_i + \sqrt{-\delta}y_i$ for $x_i,y_i \in \Q$, we have that
$q(u_1, \dots,u_d) = x_1^2+\dots+x_d^2+\delta(y_1^2+\dots+y_d^2)$, which
does not represent zero even over $\R$.
\end{proof}

\begin{proof}[Proof of Theorem \ref{IntroT3}]
Let $k'/k$ be the quadratic extension and
$q$ the anisotropic Hermitian form as in the lemma, and let
$G' = \PGU(q)$ and $G = \PGL[d]$. Then
$G'(F) \isom G(F)$,
because $k' \sub F$, so $k'\tensor[k]F = F\times F$ and
$$G'(F) = \set{(a,b)\in \GL[d](F)\times \GL[d](F)\suchthat b = a^{*}}/\mul{Z} = G(F).$$
Choose $\Gamma = G'(R_0)$ (where $R_0$ is as defined in
\Eq{R0def}) and $\valF$ is the valuation on $F$. This is a
cocompact lattice of $G'(F)$ since $G'(k) = \PGU(q,k)$ has rank
zero. Moreover if we let $q_1$ denote the sum of the first $d-1$
terms in a diagonal form of $q$, then $q_1$ does not represent
zero, and setting $H' = \PGU(q_1)$, $H'(F) = \PGU(q_1,F)$ embeds
in $G'(F)$ as $(d-1)\times (d-1)$ matrices (and is isomorphic to
$H(F) = \PGL[d-1](F)$ for the same reasons as for $q$). We deduce
that $\Lambda = \Gamma \cap H'(F)$ is a cocompact lattice in
$H'(F) = H(F)$.

By \Pref{S=sA} it remains to find a spherical non-tempered
sub-representation $\rho$ of $\L2{\dom{\Gamma_I}{G(F)}}$ for a
congruence subgroup $\Gamma_I$.

Now, since $\dom{\Lam}{H(F)}$ is compact,
$$\L2{\dom{\GL[d-1](F)}{\GL[d](F)}} = \L2{\dom{H(F)}{G(F)}} \sub \L2{\dom{\Lam}{G(F)}}.$$

The group $\GL[d](F)$ acts on $V \oplus V^*$  where $V = F^d$ and
$V^*$ is the dual space. Fixing $e_1 \in V$,
the stabilizer of $e_1 \oplus e_1^*$ is isomorphic to $\GL[d-1](F)$,
so $\L2{\dom{\PGL[d-1](F)}{\PGL[d](F)}} \sub \L2{V \oplus V^*} \isom
\L2{V \tensor F^2}$.
Now, using the action of $\GL[d]\times \GL[2]$ on $V\tensor F^2$,
one can prove that $\L2{V \tensor F^2}$ is the direct integral of
${\rho'\tensor \rho}$ over ${\mbox{tempered }\rho \in
\widehat{\GL[2]}}$ (the unitary dual), where $\rho'$ is the \rep{}
of $\GL[d]$ obtained by inducing $\rho \tensor \id_{d-2}$ from
$\GL[2]$. In particular $\rho'$ is not tempered if $d \geq 4$
(since $\id_s$ is non-tempered if $s \geq 2$). Thus,
$\L2{\dom{\PGL[d-1](F)}{\PGL[d](F)}}$ has
spherical non-tempered \subrep{s}.
We thank R.~Howe for this argument.

For an ideal $I \normali R_0$, let $\Lam_I = \Lam\Gamma(I)$. The
$\Lam_I$ have finite index in $\Gamma$ and so are cocompact in $G(F)$.
Moreover, $\cap_{I} \Lam_I = \Lam$.

Now, $\L2{\dom{\Lam}{G(F)}}$ is weakly contained in $\cup_I
\L2{\dom{\Lam_I}{G(F)}}$ \cite{BLS1},\cite{BLS2}, so for some $I \normali
R_0$, $\L2{\dom{\Lam_I}{\PGL(F)}}$ contains a spherical
non-tempered \subrep{} (which is discrete since $\Lam_I$ is
cocompact). It follows that $\dom{\Lam_I}{\B(F)}$  as well as
$\dom{\Gamma(I)}{\B}$ are non-Ramanujan.
\end{proof}

A final remark is in order: so far all the Ramanujan complexes
constructed were quotients of $\tilde{A}_{d-1}(F)$ where $F$ is an
arbitrary local field of positive characteristic. For
characteristic zero the problem is still open, except for $d = 2$.
Of course one hopes eventually to define and construct Ramanujan
complexes as quotients of the Bruhat-Tits buildings of other
simple groups as well.

\def\US{A.~Lubotzky, B.~Samuels and U.~Vishne}
\def\BIBPapI{\US, {\it Ramanujan complexes of type $\tilde{A_d}$},
Israel J. of Math., to appear.}
\def\BIBPapII{\US, {\it Explicit constructions of Ramanujan complexes},
European J. of Combinatorics, to appear.}
\def\BIBPapIII{\US, {\it Division algebras and non-commensurable isospectral manifolds},
preprint.}
\def\BIBPapIV{\US, {\it Isospectral Cayley graphs of some simple groups},
preprint.}


\begin{thebibliography}{99}

%


\bibitem[AT]{AT}
\book{E.~Artin and J.~Tate}{Class Field
Theory}{W.A.~Benjamin}{1967}

\bibitem[B1]{Ballantine}
\paper{C.M.~Ballantine}
{Ramanujan Type Buildings}{Canad. J. Math.}{{\bf 52}(6)}
{1121--1148}{2000}

\bibitem[B2]{Ballantine2}
\paper{C.M.~Ballantine}{A Hypergraph with Commuting Partial
Laplacians}{Canad. Math. Bull.}{{\bf 44}(4)}{385--397}{2001}

\bibitem[Bu]{Bump}
\book{D.~Bump}{Automorphic Forms and Representations}{Cambridge SAM {\bf 55}}
{1998}


\bibitem[BLS1]{BLS1}
M.~Burger, J.-S.~Li and P.~Sarnak,
{\it Ramanujan duals and automorphic spectrum}, unpublished,
(1990).

\bibitem[BLS2]{BLS2}
\paper{M.~Burger, J.-S.~Li and P.~Sarnak}
{Ramanujan duals and automorphic spectrum}
{Bull. AMS} {{\bf 26}(2)} {253--257}{1992}

\bibitem[C]{Cartier}
\paper{P.~Cartier}{Representations of $p$-adic Groups: A survey}
{Proc. Symposia Pure Math.}{{\bf 33}(1)}{111--155}{1979}
%

\bibitem[Cw]{Cw}
\paper{D.I.~Cartwright}
{Spherical Harmonic Analysis on Buildings of Type $\tilde{A}_n$}
{Monatsh.~Math.}{{\bf 133}}{93--109}{2001}

\bibitem[CM]{CM}
\paper{D.I.~Cartwright and W.~M{\l}otkowski}
{Harmonic Analysis for groups acting on triangle buildings}
{J.~Austral. Math. Soc. (A)}
{{\bf 56}(1)}{345--383}{1994}

\bibitem[CS]{CS2}
\paper{D.I.~Cartwright and T.~Steger}
{Elementary Symmetric Polynomials in Numbers of Modulus $1$}
{Canad. J. Math.}{{\bf 54}(2)}{239--262}{2002}

\bibitem[CSZ]{CSZ}
%
D.I.~Cartwright, P.~Sol\'{e} and A.~\.{Z}uk,
{\it Ramanujan Geometries of type $\tilde{A_n}$},
Discrete Math. {\bf 269}, 35--43, (2003).

\bibitem[Gr]{Gr}
\thesis{Y.~Greenberg}{On the spectrum of graphs and their
universal covering (Hebrew)}{Hebrew University, Jerusalem}{1995}

\bibitem[GZ]{GZ}
R.I.~Grigorchuk and A.~\.{Z}uk, \emph{On the asymptotic spectrum
of random walks on infinite families of graphs}, Random walks and
discrete potential theory (Cortona, 1997), Sympos. Math. XXXIX,
Cambridge Univ. Press, Cambridge, {188--204}, ({1999}).

\bibitem[HT]{HT}
\book{R.~Harris and R.~Taylor}{The Geometry and Cohomology of Simple
Shimura Varieties}{Anals of Math. Studies. {\bf 151}, Princeton Univ. Press}
{2001}

\bibitem[JL]{JL}
\paper{B.W.~Jordan and R.~Livne}
{The Ramanujan property for regular cubical complexes}{Duke
Mathematical Journal}{{\bf 105}(1)}{85--103}{2000}
%

\bibitem[Kn]{Knapp}
\book{A.W.~Knapp}{Representation theory of semisimple groups}
%
{Princeton Landmarks in Mathematics, Princeton University Press, Princeton, NJ}{1986}

\bibitem[L]{Lafforgue}
\paper{L.~Lafforgue}
{Chtoucas de Drinfeld et correspondance de Langlands (French)}
{Invent. Math.}{{\bf 147}(1)}{1--241}{2002}

\bibitem[LRS]{LRS}
G.~Laumon, M.~Rapoport, and U.~Stuhler,
{\it
${\mathscr D}$-elliptic sheaves and the Langlands correspondence},
Invent. Math.  {\bf 113}(2), 217--338, (1993)

%
%
%

\bibitem[Li]{Li2}
W.-C.~W.~Li, {\it Ramanujan Hypergraphs}, Geom. Funct. Anal. {\bf
14}, 380-399, (2004).

\bibitem[Lu1]{alexbook}
\book{A.~Lubotzky}
{Discrete Groups, Expanding Graphs and Invariant Measures}
{Progress in Math. {\bf 125}, Birkh\"{a}user}{1994}

\bibitem[Lu2]{Lu1}
\paper{A.~Lubotzky} {Cayley graphs: eigenvalues, expanders and
random walks} {Surveys in combinatorics (Stirling), London Math. Soc. Lecture Notes Ser.} {\bf 218}{155--189}{1995}

\bibitem[LPS]{LPS}
\paper{A.~Lubotzky, R.~Philips and P.~Sarnak}{Ramanujan graphs}{Combinatorica}{\bf 8}{261--277}{1988}

\bibitem[LSV]{paperII}
\BIBPapII

%
%

\bibitem[M]{Mac}
\book{I.G.~Macdonald}{Symmetric Functions and Hall Polynomials,
\second\ edition}{Oxford Math. Monographs}{1995}


\bibitem[Ma1]{Ma1}
\paperTrans{G.~Margulis}{Explicit group-theoretic construction of
combinatoric schemes and their applications in the construction of
expanders and concentrators (Russian)}{Problemy Peredachi
Informatsii}{{\bf 24}(1)}{51--60}
{Prob. Inform. Transmission}{{\bf 24}(1)}{39--46}{1988}

\bibitem[Ma2]{Mabook}
\book{G.~Margulis}{Discrete subgroups of semisimple Lie groups}
{Results in Mathematics and Related Areas (3), {\bf 17},
Springer-Verlag}{1991}


\bibitem[Mo]{Morg}
\paper{M.~Morgenstern}
{Existence and explicit constructions of $q+1$ regular Ramanujan graphs for every prime power $q$}
{J. Combin. Theory Ser. B}{{\bf 62}(1)}{44--62}{1994}

\bibitem[MW]{MW}
\paper{C.~M\oe{}glin and J.-L. Waldspurger}{Le spectre
r\`{e}sidual de $\GL[n]$}{Ann. Scient. \`{E}c. Norm. Sup.}{$4^{e}$
s\`{e}rie, {\bf 22}}{605--674}{1989}

\bibitem[P]{Ped}
\book{G.K.~Pedersen}{Analysis Now}{GTM {\bf 118}, Springer, New York}{1989}

\bibitem[PR]{PR}
\book{V.~Platonov  and A.~Rapinchuk}{Algebraic Groups and Number
Theory} {Pure and Applied Mathematics {\bf 139}, Academic Press}{1994}

\bibitem[Pr]{Prasad}
\paper{G.~Prasad}{Strong approximation for semi-simple groups over function fields}{Ann. Math.}{{\bf 105}}{553--572}{1977}

\bibitem[R]{Rev}
M.~Rapoport, \emph{The mathematical work of the 2002 Fields
medalists: The work of Laurent Lafforgue}, Notices of the AMS,
{\bf 50}(2), 212--214, (2003).
%

\bibitem[Ro]{Rogawski}
\paper{J.~Rogawski}{Representations of $\GL[n]$ and division algebras over a $p$-adic field}{Duke Math. J.}{{\bf 50}}{161--196}{1983}

\bibitem[Sa]{Ali}
A.~Sarveniazi,
{\it Ramunajan $(n_1,n_2,...,n_{d-1})$-regular
hypergraphs based on Bruhat-Tits Buildings of type $\tilde{A}_{d-1}$},
arxiv.org/math.NT/0401181.

\bibitem[Sc]{Schar}
\book{W.~Scharlau}{Quadratic Forms}{Queen's papers in Pure and Applied Math. {\bf 22}, Queen's Univ., Kingston, Ontario}{1969}

\bibitem[Se1]{Serre}
\book{J.-P.~Serre}{Galois Cohomology} {Springer, translated from the 1964 French text}{1996}

\bibitem[Se2]{SerreCSP}
\paper{J.-P.~Serre}{Le probl\'{e}me des groupes de congruence pour
$\SL[2]$}
{Ann. of Math.}{{\bf 92}(2)}{489--527}{1970}

\bibitem[Se3]{Serre:Trees}
\book{J.-P.~Serre}{Trees}{\second\ edition, Springer Mono. Math.}{2003}

\bibitem[T]{Tadic}
\paper{M.~Tadi\'c}{An External approach to unitary
representations}{Bull. Amer. Math. Soc.}{{\bf 28}}{215--252}{1993}

\bibitem[V]{Vig}
\paper{M.F.~Vigneras} {Correspondances entre representations
automorphes de ${\rm GL}(2)$ sur une extension quadratique de
${\rm GSp}(4)$ sur $ Q$, conjecture locale de Langlands pour ${\rm
GSp}(4)$}{Contemporary Math.}{{\bf 53}}{463--527}{1986}

\bibitem[Z]{Zelevinsky}
\paper{A.V.~Zelevinsky}{Induced Representations of Reductive $p$-adic groups II: on irreducible representations of $\GL[n]$}{Ann.~Sci.~E.N.S.}{{\bf 13}(4)}{165--210}{1980}




\end{thebibliography}
\end{document}